\documentclass{article}
\usepackage{epsfig}
\usepackage{framed}\usepackage{geometry}
\usepackage{soul,ulem}
\usepackage{amssymb}
\usepackage{amsbsy,todonotes,cancel, mathrsfs, amsmath}
\usepackage{amsthm}
\usepackage[colorlinks=true, pdfstartview=FitV, linkcolor=blue, citecolor=blue, urlcolor=blue]{hyperref}
\usetikzlibrary{calc}
\usetikzlibrary{decorations.markings}
\usepackage{verbatim}
\usepackage{ulem}    

\renewcommand{\theequation}{\arabic{section}.\arabic{equation}}
\makeatletter
\@addtoreset{equation}{section}
\makeatother

\def \scr{\mathscr}
\theoremstyle{plain}
\newtheorem{theorem}{Theorem}[section]
\newtheorem{question}[theorem]{Question}
\newtheorem{proposition}[theorem]{Proposition}
\newtheorem{corollary}[theorem]{Corollary}
\newtheorem{lemma}[theorem]{Lemma}

\theoremstyle{definition}
\newtheorem{definition}[theorem]{Definition}
\newtheorem{example}[theorem]{Example}

\newtheorem{exercise}[theorem]{Exercise}
\newtheorem{examps}[theorem]{Examples}

\newtheorem{remark}[theorem]{Remark}

\def \nn{\nonumber}
\definecolor {darkgreen}{rgb}{0,0.6,0}
\def\le{\left}
\def\eqref#1{(\ref{#1})}
\def\ri{\right}
\def\M{\mathcal M}

\def \QED{\hfill $\blacksquare$\par \vskip 5pt}
\def\ds{\displaystyle}

\def\res{\mathop{\mathrm {res}}\limits_}

\definecolor{shadecolor}{rgb}{0.95, 0.95, 0.86}

\def\br{\begin{remark}}
\def\er{\end{remark}}
\def\bt{
\definecolor{shadecolor}{rgb}{0.95, 0.95, 0.86}
\begin{shaded}
\begin{theorem}}
\def\et{\end{theorem}
\end{shaded}}
\def\bq{
\definecolor{shadecolor}{rgb}{0.95, 0.95, 0.46}
\begin{shaded}
\begin{question}}
\def\eq{\end{question}
\end{shaded}}
\def\bd{
\begin{definition}}
\def\ed{\end{definition}
}
\def \Re{\operatorname{Re}}
\def \Im{\operatorname{Im}}
\def\bp{
\begin{shaded}
\begin{proposition}}
\def\ep{\end{proposition}
\end{shaded}}

\def\bc{\begin{corollary}}
\def\ec{\end{corollary}}
\def\ord{\hbox{ ord }}
\def\brs{\begin{remarks}
\begin{enumerate}}

\def\ers{\end{enumerate}\end{remarks}}
\def\bx{\begin{example}\small}
\def\ex{\end{example}}
\def\bxr{\begin{exercise}\small}
\def\exr{\end{exercise}}
\def\bl{\begin{lemma}}
\def\el{\end{lemma}}
\def\bxs{\begin{examps}. \rm\begin{enumerate}}
\def\exs{\end{enumerate}\end{examps}}

\def\ov {\overline}

\def\&{\hspace{-15pt}&}
\def\bea#1\eea{\begin{align}#1\end{align}}
\def\beas{\begin{eqnarray*}}
\def\eeas{\end{eqnarray*}}
\def \pa{\partial}
\def\C{{\mathbb C}}

\def\R{{\mathbb R}}

\def\N{{\mathbb N}}

\def\wh{\widehat}

\def\Z{{\mathbb Z}}

\def\a{\alpha}
\def\d{\mathrm d}

\def\l{\lambda}

\def\1{{\bf 1}}

\def\wt{\widetilde}
\def\ds{\displaystyle}

\def\black{}
\def\cyan#1{{#1}\black}

\def \be#1\ee{\begin{align}#1\end{align}}

\begin{document}


\baselineskip 13pt plus 1pt minus 1pt
\vspace{0.2cm}
\begin{center}
\begin{Large}
\textbf{The Stieltjes--Fekete problem and degenerate orthogonal polynomials}
\end{Large}
\end{center}

\begin{center}
M. Bertola$^{\dagger\ddagger\star}$ \footnote{Marco.Bertola@\{concordia.ca, sissa.it\}} 
E. Chavez-Heredia$^{\ddagger\diamondsuit}$ \footnote{eduardo.chavezheredia@bristol.ac.uk}
T. Grava $^{\ddagger\diamondsuit}$ \footnote{Tamara.Grava@sissa.it}.
\\
\bigskip
\begin{minipage}{0.7\textwidth}
\begin{small}
\begin{enumerate}
\item [${\dagger}$] {\it  Department of Mathematics and
Statistics, Concordia University\\ 1455 de Maisonneuve W., Montr\'eal, Qu\'ebec,
Canada H3G 1M8} 
\item[${\ddagger}$] {\it SISSA, International School for Advanced Studies, via Bonomea 265, Trieste, Italy  and INFN sezione di Trieste }
\item[${\star}$] {\it Centre de recherches math\'ematiques,
Universit\'e de Montr\'eal\\ C.~P.~6128, succ. centre ville, Montr\'eal,
Qu\'ebec, Canada H3C 3J7}
\item[${\diamondsuit}$] {\it School of Mathematics, University of Bristol, Fry Building, Bristol,
BS8 1UG, UK}
\end{enumerate}
\end{small}
\end{minipage}
\end{center}

\begin{abstract}

A result of Stieltjes famously relates the zeroes of the classical orthogonal polynomials with the configurations of points on the line that minimize a suitable energy with logarithmic interactions under an  external field.
 The optimal configuration satisfies an algebraic set of equations:   we call this set of algebraic equations  the {\it Stieltjes--Fekete} problem. In this work we consider the Stieltjes-Fekete problem when  the derivative of the external field is  an arbitrary rational  complex  function.  We show that, under assumption of genericity,  its solutions are in one-to-one correspondence with the zeroes of certain non-hermitian orthogonal polynomials that satisfy an excess of orthogonality conditions and are thus termed ``degenerate''.  When the differential  of the external field on the Riemann sphere is of degree $3$   our result reproduces Stieltjes' original result and provides its direct generalization for higher degree after more than a century since the original result.
\end{abstract}

\vskip 2cm
\tableofcontents

%
%
\section{Introduction and results}
The weighted Fekete points   are  { the solution of the problem} defined as follows: given a smooth  (real--valued)  {\it external potential} $\mathcal Q(z,\ov z)$, $z\in\C$,  satisfying suitable additional assumptions that depend on the context, find the configuration of $n\in \N$ points  $(z_1,\dots,z_n)\in\mathbb{C}^n$ that provides 
 the maximum  of the weighted Fekete functional
\be
\mathcal F(z_1,\dots, z_n) = \prod_{j=1}^n\prod_{k=1\atop k\neq j}^n {\Big|z_j-z_k\Big|}{\rm e}^{-\frac {\mathcal{Q(}z_j) +\mathcal{Q}(z_k)}{2(n-1)}}.
\ee
 Equivalently,   the  weighted Fekete points  provide  the minimum of  the energy  functional
\be
\mathcal E(z_1,\dots, z_n) =- 2\sum_{1\leq i<j\leq n}^n\log|z_j-z_k|+\sum_{j=1}^n \mathcal{Q}(z_j)\,.
\label{holFekete}
\ee
The above energy has a  clear electrostatic interpretation: the points $(z_1,\dots,z_n)$ can be considered as charged particles in the external field $\mathcal Q$
 that   interact  {with each other } via a logarithmic potential,  (namely the force is inversely proportional to the relative distance). 
The problem of finding the critical configurations of \eqref{holFekete} is referred to as Stieltjes-Fekete problem.
An  introduction to the  weighted  Fekete problem and its connections with logarithmic potential  in  an external field can be found in \cite{SaffTotik}. 

Depending on the setup, one may require that the points belong to some assigned domain $\mathcal D$.
 A case of special interest is when the potential $\mathcal{Q}$ is a  harmonic function of the form  $\mathcal{Q} (z) =2\Re(\wh\theta(z))$ 
with $\wh\theta$ analytic, 
except for a finite number of singularities and branch cuts and with single valued derivative $\wh\theta'(z)$.

If  a configuration $\scr Z_n = \{z_1,\dots, z_n\}$  forms a {\it{critical} configuration},  it  { is then  such that } the gradient of the energy  $\mathcal E$  { vanishes}, which promptly leads to the set of equations
\be
2 \sum_{k\neq j} \frac 1{ z_j- z_k} =\pa \mathcal Q(z_j)=\wh\theta'(z_j) , \ \ \ j=1,\dots, n,
\label{criticality}
\ee
where $ \pa = \frac 1 2\le(\frac \pa{\pa x} - i \frac \pa{\pa y} \ri)$,  $z=x+i y$.
   {The following}  notable choices of $\wh\theta(z) $ are associated to  classical orthogonal polynomials:
   \begin{itemize}
   \item $\wh\theta(z) = z^2$ and $\mathcal D=\R$;
   \item  $\wh\theta(z)=  z- ( \alpha+1) \log z $,  $\alpha>-1$,  and 
  $\mathcal D=(0,\infty)$;
   \item  $\wh\theta(z)=-(\alpha+1) \log (1-z) -  (\beta+1)  \log (1+z) $,  $\alpha,\beta>-1$  and the domain  $\mathcal D=(-1,1)$.
   \end{itemize}
   {  Indeed,} for the above three choices of  $\wh\theta(z)$ and  $\mathcal D$,   the Fekete points  are  the  {zeros} of the Hermite, Laguerre and Jacobi polynomials, respectively, and they provide a global minimum to the weighted Fekete functional.
  The electrostatic interpretation of the zeros of the classical orthogonal polynomials is  due to  Stieltjes (although studied also by Bochner, Heine, Van Vleck, and Polya)  and it is 
   one of the most elegant results in the theory of special functions (for a review see \cite{Marcellanreview}).
In the present paper we choose  $\wh\theta'(z)$ to be a rational function,  so that   the condition of criticality \eqref{criticality} takes  the form
\be
\begin{split}
\label{Feketeintro} 
&{ \wh\theta'(z)=\frac{A(z)}{2B(z)}}\\
&\sum_{k\neq j} \frac 1{z_j - z_k} =\frac{A(z_j)}{2B(z_j)}, \ \ \ j=1,\dots, n,
\end{split}
\ee

where $A,B$ are two relatively prime\footnote{ { The equations clearly depend only on the ratio $A/B$; we will comment  on the assumption of being relatively prime in App. \ref{notprime}.}} polynomials  (with $B$ monic). 
Equations of similar nature as \eqref{Feketeintro} are  sometimes referred to as {\it Stieltjes--Bethe} equations because of their appearance in the Bethe-Ansatz for spin-chains \cite{Gaudin, HarnadWinternitz}   and can be considered on Riemann surfaces in higher genus \cite{KorotkinBethe}.

The question we pursue here is whether the criticality condition \eqref{Feketeintro} for \eqref{holFekete}  {bears any}  relation to zeros of orthogonal polynomials  as in  the classical case. 
In the excellent review \cite{Marcellanreview}  on this circle of ideas  the following questions were raised, to quote verbatim from loc. cit.:\\[3pt]

\noindent $\cdot$ " [...] Are there generalizations of  electrostatic models to other families of polynomials?\\
$\cdot$ Why necessarily the global minimum of the energy $\scr E$ should be considered? Which other types of equilibria described
above could be linked to the zeros of the polynomials?"\\
$\cdot$ What is the appropriate model for the complex zeros (when they exist)? [...]"
\\[10pt]
{ 
Regarding the first point, in \cite{Andrei2}  it was shown that critical configurations are related to Heine-Stieltjes polynomials, namely, to polynomial solutions of second order ODEs with polynomial coefficients. 
} 
Our present paper extends the previous literature and addresses and answers precisely the remaining of the  above three questions;   the result can be described concisely by the following statement:
\begin{shaded}
\noindent There is a   {generically}  one-to-one correspondence between the solutions of the Stieltjes--Bethe equations \eqref{Feketeintro} and
 the {\it maximally degenerate} orthogonal polynomials of degree $n$ for a {\it semiclassical moment functional} of type $(A,B)$.
\end{shaded}
 {The meaning of the qualifier ``generic'' will be spelt out in Theorem \ref{thmintro}, but we anticipate here that it is in the sense of Zariski topology (i.e. the complement of the zero  set of polynomial relations)}.

 {In the above statement}
a degenerate orthogonal polynomial is an orthogonal polynomial satisfying  an excess of orthogonality conditions  {(see Def. \ref{defdeg}).}
To further clarify the terms used above  and formulate our main theorem  (Theorem~\ref{thmintro}  below) 
we first need  to  {recall the definition of} a semiclassical moment functional \cite{Marcellan0,Marcellan,Maroni,Shohat}.
First of all, 
a  moment functional  is simply a linear map on the space of polynomials $\mathcal M:\C[z]\to \C$, associating to the monomial powers $z^j, j\geq 0$ a sequence of {\it moments} $\mu_0,\mu_1,\dots,\in \C$, and then extended by linearity to all polynomials.
\begin{definition}
\label{defsemi}
A moment functional $\mathcal M:\C[z]\to \C$  is called {\it semiclassical} if there exist two  polynomials 
$A(z), B(z)$ of degree $a,b$, respectively such that 
\be
{\mathcal M}\big[A(z)p(z)\big] ={ \mathcal  M}\big[B(z) p'(z)\big],\quad \ \ \forall p(z) \in \C[z].
\ee
We say  that such a semiclassical functional $\mathcal M$ is of type $(A, B)$.
\end{definition}
 {In our application the polynomials $A,B$ are relatively prime given the nature of the equations \eqref{Feketeintro}, but we will comment in Appendix \ref{notprime} on the nature of this assumption and how to (generically) lift it.}
\noindent The  concept of {semiclassical}  moment functional  originated in \cite{Maroni0, Maroni},  and it was then developed \cite{Marcellan}, see also the review \cite{Marcellanreview}.

The main result  of \cite{Ismail, Marcellan,Maroni} (see also the introduction of \cite{Bertosemiclassical} for a quick rundown) is that  a  semiclassical moment  functional  $\mathcal M$
can be represented by the  moment functional  $\mathcal M_{\Gamma,\theta}$  defined  through the expression
\be
\label{moment}
\mathcal M_{\Gamma,\theta}:\C[z]&\to \C,
\nn \\
z^j&\mapsto \mathcal M_{\Gamma,\theta}[z^j]= \mu_j =\int_{\Gamma} z^j{\rm e}^{\theta(z)} \d z,
\ee
where the  {\it symbol}\footnote{ {The term ``symbol'' is used here in parallel with the use in the theory of T\"oplitz matrices where the same moments (allowing also the negative ones) are arranged in a namesake matrix rather than in a matrix of Hankel type}.}  $\theta$ of the exponential weight  satisfies
\be
\theta' (z) =- \frac {A(z) + B'(z)}{B(z)}\,,
\label{semidef1}
\ee
and the integral is taken over a ``weighted contour'' $\Gamma=\sum_{j=1}^ds_j\gamma_j$ in the following sense:
\begin{equation}
    \int_{\Gamma}f(z) \d z := \sum_{j=1}^d s_j \int_{\gamma_j} f(z) \d z.
\end{equation} 
The weighted contour $\Gamma=\sum_{j=1}^ds_j\gamma_j$ is expressed in terms  of  contours $\gamma_j$ in the complex plane that extend from a zero of $B$ to another (or to infinity) described in Section \ref{contourssec}; the complex parameters $s_j, \ j=1,\dots d$ parametrize the space of semiclassical moment functionals of a given type $(A,B)$.
 Examples of (semi)classical moment functionals are in Table \ref{classical}.
\begin{table}[t]
\centering
\begin{equation*}
\begin{array}{c|c|c|c}
\text {Name} & \text{Symbol } \theta(x) & \text{Type } (A,B) & \text{Contour } \Gamma \\[5pt]
\hline & & & \\
\text {Jacobi} & \alpha \log (1-x) +  \beta  \log (1+x)  &  -\left(\beta +\alpha +2) x -(\alpha-\beta), x^2-1\right) &   [-1,1] \\[10pt]
&\alpha,\beta>-1&&\\
\hline & & & \\
\text{Hermite} &  - {x^2} &  (-2x,1) & \R \\[10pt]
\hline & & & \\
\text{Laguerre} &  - x+ \ \alpha  \log x,\;\;\alpha>-1  &  (x-\alpha -1, x) & \R_+ \\[10pt]
\hline &&&\\
\color{black} \text{Bessel} &-\frac 1 x + \nu \ln x,\;\nu\in\C &( 1+(2- \nu) x, x^2) &\gamma_1 \color{black}
\end{array}
\end{equation*}
\caption{
The classical orthogonal polynomials can be obtained from a semi-classical moment functional of type $(A,B)$ with symbol $\theta(x)$ and contour $\Gamma$ as indicated.
The common feature is that $\d \theta = \theta' \d x$ is a differential on the Riemann sphere $\mathbb P^1$ with total degree of poles $=3$. 
 The last line is the class considered in \cite{KrallBessel}: in this case  the contour $\gamma_1$ can be chosen as a cardioid with the cusp at the origin and the lobe on the right, or a circle centered at the origin if $\nu\in \Z$. 
}
\label{classical}
\end{table}
The maximum  number $d$ of  { linearly independent}  contours is the degree of the pole divisor of $\theta'(z) \d z$ on the Riemann sphere minus $2$.
In the following, to keep notation simple, we omit the explicit indication  of dependence of the moment functional  on $\Gamma$ and $\theta$ and write simply $ \mathcal M$ in place of $\mathcal M_{\Gamma,\theta}$.

Given a   semiclassical moment functional $\mathcal M$ the corresponding  orthogonal polynomials, when they exist, are a sequence $\{P_j(z)\}_{j\in \N}$ of polynomials, each of degree $\leq j$, satisfying
\be
\label{orthointro}
\mathcal M\big[P_j(z)P_k(z)\big]= \delta_{jk} h_k, \quad  j,k\in \N,
\ee
for some numbers  $h_k\in \C {\setminus \{0\}}$.
We  { now need to }  introduce the concept of $\ell$-degenerate orthogonality. \bd
\label{defdeg}
A polynomial $P_n$ of degree $n$  is called {\bf $\ell$--degenerate orthogonal}, $\ell\geq 0$,  if it satisfies the following excess of orthogonality conditions
\begin{equation}
\label{deg_eq}
{ \mathcal M}\big[P_n(z)z^{k}\big]= \sum_{j=1}^d s_j \int_{\gamma_j} P_n(z)z^{k}{\rm e}^{\theta(z)} \d z=0, \ \ \ k=0,1\dots, n+\ell-1.
\end{equation}
The polynomial $P_n(z)$, is called maximally degenerate if $\ell=d-1$ with  $d=\max\{\deg A, \deg B-1\}.$
\ed
\cyan{For a chosen value of $n$ the condition that $P_n$ is a  maximal degenerate  polynomial, requires a specific choice of  the weights $s_1,\dots,s_d$.   Once $P_n$
is chosen  to be  maximal degenerate,   then $P_m$ for $m\neq n$ is not   in general a  maximal degenerate orthogonal polynomial.  To stress the dependence on $n$ of the moment functional $\mathcal M$  we denote it  by $\mathcal M_n$. }

 {
We recall  \cite{AptMarRo} that the integer $d-1$ is referred to as the ``class'' of the semiclassical moment functional, and hence we may say that a maximally degenerate orthogonal polynomial satisfies a number of additional orthogonality relations equal to the class of the moment functional.}
 For a given type of moment functional
   and given degree $n\in \N$, we consider the conditions of degeneracy as (homogeneous) constraints on the parameters $s_1,\dots, s_d$. For this reason in general we expect a maximal degeneracy of $d-1$: see Lemma \ref{lemmaell}. 
 If  $d=1$ then  any orthogonal polynomial is maximally degenerate (i.e. $0$--degenerate) by default.   This applies to 
all the classical moment functionals  (see Table \ref{classical}).
With these preparations, we can formulate  our main results.
\bt
 \label{thmintro}
 {    Let $A$ and $B$ be relatively prime  polynomials with    $B$ having  at least one multiple root  if $\deg B<\deg A$. Further let  $d=\max\{\deg A, \deg B-1\}$ and 

  \begin{equation}
  \theta' (z) =- \frac {A(z) + B'(z)}{B(z)}\,. 
  \end{equation}
  Suppose that for given integer $n$ there exist constant parameters $s_1,\dots, s_d$ such that  the  semiclassical  moment functional $\mathcal M$ of type $(A,B)$    admits a maximally degenerate polynomial $P_n$ of degree $n$  according to  definition~\ref{defdeg}.}
   Then 
 \begin{enumerate}
 \item [{\bf (1)}] the function $F(z) = \sqrt{B(z)} P_n(z) {\rm e}^{\frac{1}{2}\theta(z)}
 $ solves the differential equation 
 \be
 F''(z) - V(z) F(z) =0\,,
 \ee
 where the potential $V(z)$ is a rational function with poles only at the zeros of $B(z)$ at most of twice the order and takes  the form:
 \be
 V(z) =  
 \frac 1 2\left(\theta' +\frac {B'}{B}\right)'+ \frac 1 4\left(\theta'+\frac {B'}{B}\right)^2 
  + \frac {Q}{B} , \ \ \ \deg Q \leq d-1,
 \ee
 with $Q$ a polynomials of $\deg Q \leq d-1$.
 Equivalently the polynomial $P_n$  satisfies  the  generalized Heine-Stieltjes  equation
 \be
 \label{ODEPn0}
  B(z) P_n'' - A(z)  P_n' - Q(z) P_n(z) = 0 .
 \ee
 \item [{\bf (2)}] Let  $(z_1, \dots, z_m)$ denote the roots of $P_n$  { that are   distinct  from those of $B$};  they  satisfy the system of equations
 \be
 \sum_{\ell\neq j\atop 1\leq \ell\leq m }\frac 1{ z_j-z_\ell}   = \frac {A(z_j) }{2B(z_j)}
 - \sum_{{ c: \atop B(c)=0=P_n(c)}} \frac {\sharp_c}{z_j - c},\ \ \ j= 1,\dots, m\,,
 \label{Fekete0}
 \ee
 { and  the positive integer $\sharp_c$ denotes the multiplicity of the root $c$ of $B(z)$.}
  {
 Under the genericity assumption that all roots  $c$ of $B$  satisfy the condition $A(c) - kB'(c)\neq 0$ for  all $k =0,1, \dots$, the polynomial $P_n$ does not share any root with $B$}  {  and all its zeros $\{z_1,\dots,z_n\}$ satisfy the system of equations
 \be
\label{Fekthm}
\sum_{\ell\neq j\atop 1\leq \ell\leq n} \frac 1{z_j - z_\ell} = \frac{A(z_j)}{2B(z_j)}, \ \ \ j=1,\dots, n.	
\ee
Viceversa let $\scr Z= \{z_1,\dots, z_n\}\in\C^n$ be a critical configuration satisfying the equilibrium equations \eqref{Fekthm}, 
  then  the polynomial $P_n(z) = \prod_{j=1}^n (z-z_j)$ is a (non-hermitian) maximally degenerate orthogonal polynomial  for a  semiclassical moment functional $\mathcal M  {=\mathcal M_n}$  of type $(A,B)$ and item  {\bf (1)}  holds. In the case $\deg A<\deg B$ the additional assumption on $n$ has to be imposed
  \[
  2n > \Re \Lambda -1 - \deg B, \ \ \ \ \Lambda:= - \res{z=\infty} \theta'(z)\d z.
  \]
}
 \end{enumerate}
 \et
\noindent The proof of the above theorem is contained in Theorems \ref{thm1} and  \ref{thm2}.

For the classical case ($d=1$)  maximal degeneracy is the $0-$degeneracy  and  the statement of our theorem  overlaps entirely with the classical result \cite{Stieltjes,Bochner} (see Table~\ref{classical}).
The differential equation that appears in Theorem \ref{thm1} is simply the (stationary) Schr\"odinger equation  form of the second order ODE satisfied by the classical orthogonal polynomials, which also characterizes them as shown by Bochner in \cite{Bochner}, a result later generalized by Krall-Littlejohn \cite{Krall}. 
When the parameters $\alpha,\beta\leq -1$ in the Jacobi or Laguerre weight,  the corresponding Stieltjes-Bethe equation has appeared in the study of the  weakly
anisotropic Heisenberg spin chain \cite{SD, YM}. The solution of the Stieltjes-Bethe is given by the zeros of the generalized Jacobi and Laguerre polynomials. 
In this latter case  the asymptotic location of zeros  in the complex plane has been studied in \cite{KMc} and \cite{add1,add2,add3,add4}.
{ 
\begin{remark}
\label{remarkprime}
While the equations \eqref{Fekthm} are  invariant under simultaneous multiplication of $A$ and $B$ by {  a  polynomial $q$},  on the other hand the expression of $\theta'$ and hence the orthogonality  conditions on $P_n$ are not invariant. This observation notwithstanding,  we show in Appendix \ref{notprime} that (generically) the {\it same} polynomial $P_n$ is orthogonal for the new moment functional. Namely, the theorem still generically holds if we lift the (simplifying) requirement that $A,B$ are relatively prime.  Here ``genericity" means that none of the points $z_1,\dots, z_n$ is a zero {  of the polynomial  $q$}  by which we multiply $A$ and $B$. See App. \ref{notprime} for more detail.
\end{remark}
}

\paragraph{Example of semiclassical  functional: Freud weight  $\theta(z)= -\frac{z^4}2$ \cite{Freud}.}

We   consider the weighted contour $\Gamma=s_1\gamma_1+s_2\gamma_2+s_3\gamma_3$  where $\gamma_j$ is the contour extending to infinity from  the  directions $\arg(z)= (j-1)\frac {\pi}2$ to  $\arg(z)= \frac {j\pi}2$,   $j=1,2,3$,  and the  complex  parameters $s_\ell$ are not all simultaneously zero.  Next we define 
the  moment functional \cite{Chihara,SzegoBook} 
\be
\mathcal M[z^j] = \le(s_1\int_{\gamma_1}  + s_2\int_{\gamma_2} + s_3\int_{\gamma_3} \ri)z^j{\rm e}^{-z^4/2}\d z\,.
\label{Freudint}
\ee
A characterization \cite{Marcellan} of such a  moment functional is that it satisfies the semiclassical condition
\be
\mathcal M[2z^3 \,  p(z) ] = \mathcal M[p'(z)],\ \ \ \forall p\in \C[z]\label{Freud4}
\ee
that is a simple consequence of integration by parts, so that $\mathcal M$ is  a semiclassical moment functional of type $(A(z)=2z^3, B(z)=1)$. The parameters $s_\ell$ can be thought of as parametrizing the space of solutions of \eqref{Freud4}.

The   corresponding orthogonal polynomials, when they exist, are a sequence $\{P_j(z)\}_{j\in \N}$ of polynomials of degree $\leq j$ that satisfy  the orthogonality relation 
\eqref{orthointro}.
Clearly only the ratios of the parameters $s_\ell$ are relevant for the definition  of orthogonal polynomials,  so that we can think of the family of functionals \eqref{Freudint} with the  same symbol  $\theta(z)=-\frac{z^4}{2}$  as parametrized by $[s_1:s_2:s_3]\in\mathbb P^2$.    Having the freedom to choose the point $[s_1:s_2:s_3]\in\mathbb P^2$, one can impose 
an excess of orthogonality.
In this case the notion of {\it maximally degenerate} orthogonal polynomial  is the following;   for any  $n\in\N$   there is $[s_1:s_2:s_3]\in\mathbb P^2$  { and } a (nontrivial) polynomial $P_n(z)$ of degree $ n$ with the properties 
\be
\mathcal M[z^j P_n(z) ]= 0,\ \ \ j=0,1,2,\dots, n-1,{\bf n, n+1}.
\ee
We have emphasized with the bold font  that the orthogonality of $P_n$ extends beyond the range of powers that characterizes an ordinary orthogonal polynomial. The reader with some prior experience  {may anticipate} that the two extra conditions can be fulfilled if and only if  two certain determinants of size $n+1$ that involve the moments vanish simultaneously: this is indeed the case, see Lemma \ref{lemmaell}.   This places two homogeneous polynomial conditions on the parameters $s_1,s_2,s_3$. 
In this case our theorem states that the zeros of such a maximally degenerate polynomial will satisfy eq. \eqref{Feketeintro} with $A(z) = 2z^3= -\theta'(z)$ and  $B=1$.  

\begin{figure}
\begin{center}
\includegraphics{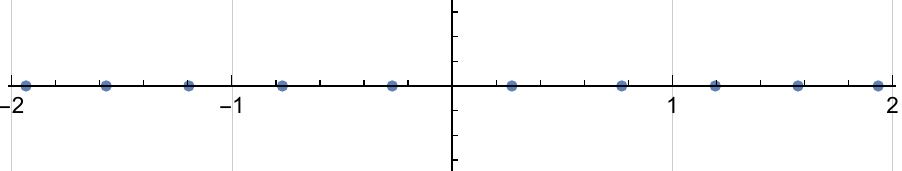}
\includegraphics{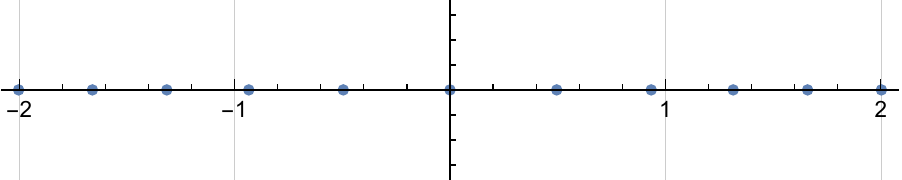}
\end{center}
\caption{The numerically computed Fekete points with $n=10, 11$ on the real axis for the Freud weight ${\rm e}^{-x^4/2}$.}
\label{Feketeintroreal}
\end{figure}

\begin{remark}
\label{rm_zeros}
\rm
{The   function $\wh \theta(z)=\int \frac {A}{B}$ appearing in the  energy functional \eqref{holFekete}  is just $\wh \theta(z) = z^4/2=-\theta(z)$,  and it is real for $z\in \R$; }it is simple to see from the electrostatic interpretation that there is an optimal configuration  {minimizing the Fekete functional} with  $z_j\in\R$.
Our theorem says that the corresponding polynomial $P_n(z) = \prod_{j=1}^n(z-z_j)$ is indeed an orthogonal polynomial, but {\it not} for the 
 orthogonality  chosen on the real axis  (this could not be because then the moment functional is strictly positive definite and the Hankel determinant of the moments cannot vanish).  
 The moment functional is 
\be
\begin{split}
&\mathcal M[z^{2j+1}] = \int_\R z^{2j+1}{\rm e}^{-\frac {z^4}2}\d z + s \int_{\i \R} z^{2j+1}{\rm e}^{-\frac {z^4}2}\d z=0\\
&\mathcal M[z^{2j}] = \int_\R z^{2j}{\rm e}^{-\frac {z^4}2}\d z + s \int_{\i \R} z^{2j}{\rm e}^{-\frac {z^4}2}\d z= (1- (-1)^j \Im (s)) 2^{\frac{2j-3}{2}}\Gamma\left(\frac{2j+1}{4}\right)\\
\end{split}
\ee
 where  $\Gamma$ is the standard gamma-function, and $s$ has been chosen pure imaginary in order to obtain real moments.
 We observe that, using our previous notation,   the integration over the real line $\R$ is homotopic to the path $-\gamma_1-\gamma_2$, while the integration along 
 the imaginary axis is homotopic to the path  $-\gamma_2-\gamma_3$.
 The degeneracy condition on the polynomial of degree $n=10$, gives $s=0.00001349595 \,i $ and for $n=11$ one obtains $s=-3.79352745 \,10^{-6}\,i $.
 The zeros of the corresponding degenerate orthogonal polynomials 
give a  solution of the   Stieltjes--Bethe equation \eqref{Feketeintro}, and are plotted   in Fig. \ref{Feketeintroreal}.
\hfill $\triangle$
\end{remark}
 \begin{remark}
 We remark that for real  continuous functions  $\mathcal{Q}:\R\to\R$ such that 
$ \lim_{|x|\to\infty}\frac{\mathcal{Q}(x)}{\log |x|}=+\infty$,    the weighted Fekete problem  on the real line
\[
\mathcal E(x_1,\dots, x_n) =- 2\sum_{1\leq i<j\leq n}^n\log|x_j-x_k|+n\sum_{j=1}^n \mathcal{Q}(x_j)\,,
\]
 is solved, in the limit $n\to\infty$, by the minimizer $\mu^*$  of the  convex  functional 
\[
\int\int\log\frac{1}{|s-t|}d\mu(s)d\mu(t)+\int\mathcal{Q}(s)d\mu(s)
\]
with respect to  probability measures $d\mu$  on the real line. Such a  problem has been extensively  studied and it is known that  the zeros of the orthogonal 
polynomial $p_n(x)$  with respect to the exponential weight  $e^{-n\mathcal{Q}(x)}dx$  converge in the limit $n\to\infty$ to the same measure $\mu^*$, (see e.g. \cite{MhaskarSaff, GoncharRakhmanov, Deift1, Deift2, KMcL, SaffTotik}).\hfill $\triangle$
 \end{remark}

\paragraph{Connection with recent results.}
  H.~E.~Heine \cite{He} and T.~J.~Stieltjes \cite{Stieltjes}   consider the electrostatic problem with  $p+1$ fixed charges on the real line at position $b_0,\dots, b_p$ with  positive  mass $\nu_0/2,\dots \nu_p/2$  and $n$ movable charges. In this case the energy takes the form 
\begin{equation}
\label{HS}
 {\mathcal E}(z_{1}, \ldots, z_{n}) = -2 \sum_{1\leq j <k \leq n}  \log
{|z_{j}-z_{k}|} - \sum_{j=1}^{n} \sum_{\ell=0}^{p} \nu_\ell  \log
{|z_{j}-b_\ell|}\,.
\end{equation}
The problem of finding the critical points of the above energy is now commonly known as 
the  {\it classical Heine-Stieltjes electrostatic problem}. Besides Heine and Stieltjes,  the problem  has been considered by   E.~B.~Van Vleck \cite{VVl},  G.~Szeg\"{o} \cite{SzegoBook},  and G.~P\'olya \cite{Pol}, just to mention a few of the classical references. The case of negative masses $\nu_j$ has been  recently considered  in \cite{DV}.
 The  Fekete  energy  \eqref{HS} and its generalizations called  ``master functions'' were thoroughly studied by A. Varchenko and coauthors in a large number of publications including  \cite{Varchenko}.   A recent survey  can be found in the excellent review by B. Shapiro \cite{Boris11_1}. 
In the recent paper \cite{Boris18} the authors consider  {the} problem of finding the critical configurations of the Heine-Stieltjes functional \eqref{HS}
{\it subject to constraint} $$\sum_{j=1}^n r(z_j)=0,$$ where $r'(z)$ is a rational function with $B(z):= \prod_{\ell=0}^p (z-b_\ell)$ in the denominator.  {\color{black}Here $b_\ell\neq b_j$ 
for $\ell\neq j$. } This polynomial $B$ plays the same role as our $B$.\footnote{In loc. cit. the use of $A, B$ is swapped but for ease of the discussion we present their result here in a way that better overlaps with our notation.}
This problem is very close to ours; by the method of Lagrange multiplier, their  equation 
 expressing the stationarity can be cast as 
 \be
 \sum_{k\neq i} \frac 1{z_i-z_k} =  - \sum_{\ell=0}^{p} \frac {\nu_\ell}{z_i -b_\ell}  + \rho^* r'(z_i),
 \label{lagrange}
 \ee
 for an appropriate Lagrange multiplier $\rho^\star$.  They find that the critical configurations are in one--to--one correspondence with a polynomial solution $y(z) = \prod_{j=1}^n(z-z_j)$ of a suitable  {\it degenerate} Lam\'e\ equations (in their terminology) 
 $$
 B(z)y'' + \wt A(z)y' + Q(z) y=0, \ \ \ \ \wt A(z) = -\rho^\star B(z) r'(z) + A(z)
 $$
 where $A(z)$ is the numerator in the fraction $\frac {A(z)}{B(z)} = \sum_{\ell=0}^p \frac {\nu_\ell} {z-b_\ell}$. The function  $Q(z)$ is an appropriate (implicitly defined) polynomial of degree at most $\max\{\deg B-2, \deg \wt A-1\}$; this should be compared with the expression in \eqref{ODEPn}. 
 We can see that the equation \eqref{lagrange} can be recast into an  equation similar to \eqref{Feketeintro}, with the polynomial case corresponding to the case $p=0$.  The difference in our approach is that {\bf (i)} we do not have a constrained equilibrium problem;  {\bf (ii)} we can have poles of arbitrary order (i.e. the multiplicities of the zeros of $B$ are arbitrary);
 {\bf(iii)} we  connect the solution of the Stieltjes-Bethe equation with (degenerate) orthogonal polynomials, that is the main point of our work and that generalizes Stieltjes  seminal  result.

 Another complementarily related result is the work \cite{Andrei} where, for an {\it arbitrary} polynomial $P_ n(z)$ and  semiclassical weight $w$ they show, amongst other results, how to construct an ``electrostatic partner'' polynomial $S$ of degree $n-1$ so that the zeros of $P_n$ satisfy an electrostatic equilibrium problem similar to \eqref{Feketeintro} where on the right side, in addition to the logarithmic derivative of the weight, there is the contribution of $S'/S$. 
{   The electrostatic partner $S$  looks like the Wronkstian determinant defined in proposition~\ref{Wronsk}  so that when the polynomial $P_n$  is maximally degenerate, the electrostatic partner becomes constant,  and the results in \cite{Andrei}  match with our results.}

 Finally another result that connects orthogonal polynomials to potential theory  but  {goes in the opposite direction relative to } our result  is due to M. Ismail \cite{Ismail2020}:
let $p_n$ be orthogonal polynomials on $[c_0,c_1]\subseteq \R$  with respect to   an exponential weight $e^{-v(x)}$.
Then the zeros of $p_n(x)$ have an electrostatic interpretation as the critical points of the energy functional 
$\scr E(x_1,\dots, x_n) =-2 \sum_{1\leq i<j\leq n}^n\log|x_j-x_k|-\sum_{j=1}^nQ(x_j)$ where  the external potential $Q(x)=-v(x)+\log(A_n(x)/a_n)$ and $A_n(x)$ is obtained from $p_n(x)$ as explained in that paper.

\section{Semiclassical moment functionals}
\label{semiclasssec}
We now provide a quick summary of the definition and some properties  of semiclassical moment functionals. The notion originates in the work of  {Shohat \cite{Shohat} and then developed by } Maroni \cite{Maroni} and Marcell\'an-Rocha \cite{Marcellan}, and has also been extended to the bi-orthogonal case in \cite{Bertosemiclassical}. 

{\color{black} Recall now that for  two  polynomials $A(z),  B(z)$ of degree $a,b$, respectively,
the  {\it semiclassical moment functional} of type $(A,B)$ are those satisfying  Definition \ref{defsemi}.}
The main result  of \cite{Ismail, Marcellan,Maroni} (see also the introduction of \cite{Bertosemiclassical}) is that any such moment functional can be represented in a similar form as \eqref{Freudint}: 
\be
\mathcal M[p] = \sum_{\ell=1}^d s_\ell \int_{\gamma_\ell} p(z){\rm e}^{\theta(z)} \d z\ ,\qquad 
\theta' (z) =- \frac {A(z) + B'(z)}{B(z)}
\label{semidef}
\ee
and $\gamma_\ell$ are suitable contours that extend from a zero of $B$ to another (or to infinity) and described in the next section. 
If follows from a judicious use of Cauchy integral theorem and counting that there are $d = \max\{b-1,a\}$ such linearly independent contours; this integer $d$ is also the total degree of the poles of the meromorphic differential  $\theta'(z) \d z$  on the Riemann sphere, minus $2$. 

For brevity we will denote by $\Gamma = \sum_{\ell=1}^d s_\ell \gamma_\ell$ the element of a suitable homology space and simply denote by $\int_\Gamma$ the operation $\sum s_\ell \int_{\gamma_\ell}$. 
\begin{remark}
\label{catch22}
We warn the reader that the symbol $\theta$ as defined in \eqref{semidef} may be regular at some (or all) of the zeros of $B$ if they simplify in the ratio that appears in \eqref{semidef}; in particular this means that there may be different moment functional with the same symbols but a different number $d$ of contours of integration. This happens when  there are ``hard-edge'' contours of integration.  An  example  is the case of the Jacobi polynomials with $\alpha=\beta=0$; in this case the symbol is $\theta=0$ and the moment functional $\mathcal M [p] = \int_{-1}^1 p(x) \d x$ satisfies
\be
\mathcal M\Big[(x^2-1) p'(x) \Big]=  \mathcal M\Big[-2x p(x)\Big] ,\ \ \ \forall p\in \C[x],
\ee
so that $B= x^2-1$ and $A=-2x$. This follows from a  simple integration by parts, where the contribution of the boundary evaluation vanishes thanks to the fact that $B(\pm1)=0$. 
\end{remark}

\subsection{The contours of support of $\Gamma$: homology and intersection pairing}
\label{contourssec}
To provide a complete and general description of the moment functionals we need to introduce a notion used in \cite{Bertosemiclassical} which defines  dual homologies and a non-degenerate intersection pairing. The notion is, interestingly, related to the notion of {\it bilinear concomitant} of a pair of differential equations, one the (formal) adjoint of the other; details of this connection are in \cite{Bertosemiclassical}.

Before addressing the full generality of the issue  let us illustrate the construction for  the simple but non-trivial example of weights of Freud type  {  (see also \cite{AptMarRo} for further reference  on the construction of the contours $\gamma_j$)}.
\paragraph{Guide example: Freud--type weights.}
Consider the case where $B(z)=1$ and $A(z)$ is a polynomial of degree $a$, so that $\theta(z) =- \int A \d z$ is a polynomial of degree $a+1$ and $d = a$. Without major loss of generality we assume that $A(z)= z^a + \dots$ is a monic polynomial. 

\begin{figure}
\begin{center}
\resizebox{0.55\textwidth}{!}{
\begin{tikzpicture}[scale=0.9]
\draw [dashed] (0,0) to (0:6);
\draw [dashed, blue] (0,0) to (36:6);
\foreach\a in{0,1,2,3,4,5,6,7,8,9}
{\node at (\a*36:6.5) {$\infty^{(\a)}$};}

\foreach \a in{1,2,3,4}{
\draw [dashed] (0,0) to (\a*72:6);
\draw [dashed, blue] (0,0) to (\a*72+36:6);
\draw [line width=1, postaction={decorate,decoration={markings,mark=at position 0.3 with {\arrow[line width=1.5pt]{>}}}}] (\a*72-72+2:6) to  [out=180+\a*72-72, in= 180+\a*72]node[pos=0.3,above, sloped] {$\gamma_\a$}  (\a*72 -2:6);

\draw [blue,line width=1, postaction={decorate,decoration={markings,mark=at position 0.6 with {\arrow[line width=1.5pt]{>}}}}] (36-72+5-2*\a:6) to  [out=180-72+36, in= 180+\a*72-36 +2*\a-5  ]node[pos=0.7,above, sloped] {$\wh \gamma_\a$}  (\a*72 -36+2*\a-5:6);
} 
\fill [fill=white] (0,0) circle[radius=12pt];
\end{tikzpicture}}
\end{center}
\caption{The contours, and dual contours, for the Freud case, $\theta(z)= -z^5$, here with $d= a = 4$. In the figure  $\infty^{(j)}$ are
 the asymptotic directions $\arg(z) = \frac {j \pi}{5}$ for $j=0,\dots, 9$.  }
\label{Freudfig}
\end{figure}
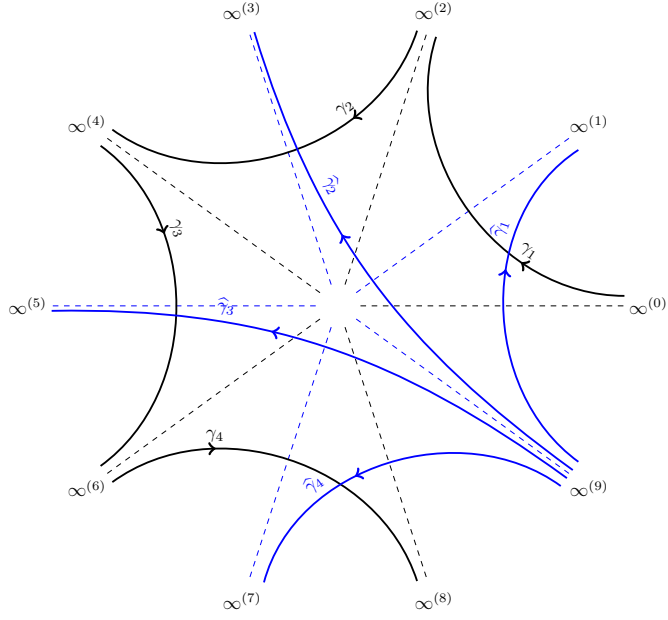

The linear space of moment functionals is in one-to-one correspondence with the solutions of the Pearcey--like ODE
\be \label{pearcey}
A(\pa_\l) \Phi(\l) = \l \Phi(\l).
\ee
Its solutions are  described as follows. We denote by $\infty^{(j)}$ the asymptotic directions $\arg(z) = \frac {j \pi}{a+1}$ for $j=0,\dots, 2a+1$ and by  {\color{black}  $\gamma_{j+1}$ the oriented contour connecting $\infty^{(2j)}$  to $\infty^{(2j+2)}$ for $j=0,\dots a$ and by $\gamma_{a+1}$  the oriented contour from $\infty^{(2a)}$ to $\infty^{(0)}$. } Then the space of solutions to the ODE \eqref{pearcey} is spanned by
\begin{equation}
    \Phi_j(\l) = \int_{\infty^{(2j)}}^{\infty^{(2j+2)}} \hspace{-10pt} {\rm e}^{\theta(z) + \l z} \d z,\quad  j=1,\dots ,a\,.
\end{equation}
Furthermore we associate solutions of the Pearcey equation to moment functionals by 
\begin{equation}
\Phi_j(\l) \ \ \longleftrightarrow \ \ \ \mathcal M_{\gamma_j}[p] = \int\limits_{\gamma_j}  p(z){\rm e}^{\theta(z)} \d z,\ \ p\in \C[z],
\end{equation}
where throughout, the contour integrals are  seen to be absolutely convergent on account of the fact that $\theta = -\frac {z^{a+1}}{(a+1)} + \mathcal O(z^a)$.
Due to the Cauchy residue theorem the $a+1$ moment functionals defined in the above formula satisfy the linear relation 
\be
\mathcal M_{\gamma_1}+ \dots + \mathcal M_{\gamma_{a+1}} \equiv 0
\ee
and hence the general moment functional of type $(A,1)$ is expressed as 
\be
\mathcal M[p] =\sum_{j=1}^a s_j  \mathcal M_{\gamma_j}[p]=\sum_{j=1}^a s_j \int_{\gamma_j}  p(z){\rm e}^{\theta(z)} \d z,\ \ \forall p\in \C[z].
\ee
The {\it dual contours} $\wh \gamma_j$ are the contours extending from $\infty^{(2a+1)}$ to $\infty^{(2j-1)}$, $j=1,\dots, a$, see Fig. \ref{Freudfig}.
They have the property that the intersection number is 
\be
{\color{black} \wh  \gamma_j\circ \gamma_\ell = \delta_{j\ell}.}
\ee

\cyan{   We recall that given two oriented curves $\gamma$ and $\eta$, their intersection number  $\gamma\circ\eta$ counts  the number of points of intersection, each counted with a $+1$ if the tangent vectors of the curves $\gamma$ and $\eta$  at the intersection point form a positively oriented frame, and $-1$ otherwise.}

\subsubsection{Case when $\deg A\geq \deg B$  with $A$ and $B$ relatively prime}
We report a canonical choice of contours $\gamma_1,\dots, \gamma_d$ as defined in \cite{Bertola:Semiiso, Bertosemiclassical, MillerShapiro}.
We say that a zero of $B$ is a {\it visible singularity} if it is a pole of $\theta'$; the zeros of $B$ that simplify in the ratio that defines  the symbol \eqref{semidef} will be called {\it hard-edges}. Note that these latter are necessarily simple zeros of $B$ due to the assumption that $A,B$ are relatively prime.
\be
\theta' (z) = -\frac {A(z)+B'(z)}{B(z)}= U_\infty'(z) + \sum_{j=1}^k U_{j}'\le(z\ri)
\ee
where $U_\infty'(z)$ is a polynomial of degree $d_\infty-1$ with  $d_\infty:=\deg A-\deg B+1$,   and 
$U_j'(z)$ are polynomials in $(z-c_j)^{-1}$ of degree  at most $ \ord_{c_j} B$, where $c_j$ are the visible singularities.  If $\deg A<\deg B$, then the term $U_\infty'(z)$ is simply zero.
   The local behaviour of $\theta(z)$ for each visible singularity $c=c_j$ of order $d_{c_j} {:= \ord_{c_j} B-1\geq 0}$ is 
\be
\theta(z) = \frac {T_j}{(z-c_j)^{d_{c_j}} }(1 + \mathcal O(z-c_j)) + r_j \ln (z-c_j), \quad z \to c_j,
\ee
 where 
\[
r_j=\res{z=c_j} \theta'(z) \d z.
\]

\begin{definition}[Local directions of steepest descent] \label{directions}
Let $c \in \mathbb{P}^1$ be a pole of $\theta'(z)\d z$.
\begin{itemize}
    \item 
    Suppose that $c$ has order $d_c+1 \ge 2$ so that 
    \begin{equation}
        \theta'(z)\d z = \frac {T_c}{\zeta_c^{d_c+1}}(1 + \mathcal O(\zeta_c) )\d \zeta_c
    \end{equation}
       where $T_c$ is the coefficient of the leading singularity in the local parameter $\zeta_c$.
   
    \footnote{To avoid lengthy case distinctions, the local parameter near a point $c\in \C$ is simply 
    $\zeta_c = z-c$ while if $c=\infty$ the local parameter will be $\zeta_\infty  = \frac 1 z$.}
       We denote by 
    \begin{equation}
    c^{(\ell)} := \le\{\arg(\zeta_c) = \frac{ - \arg T_c}{d_c}  + \ell \frac {\pi} {d_c}\ri\}, \  \ \ell = 0,1,\dots,   2d_c-1.     \end{equation}
    Namely, $\Re \theta$ tends to $-\infty$ along the directions $c^{(2k)}$ and to $\infty$ along the odd directions $c^{(2k+1)}$. These will be called {\rm local directions of steepest descent, ascent} (respectively). Note that $c^{(2d_c)}=c^{(0)}$.

    \item 
    Suppose that $c$ is a simple pole. We will denote its residue by 
    \begin{equation}
        r_c= \res{c} \theta'(z) \d z.
    \end{equation}
    A simple pole with $\Re (r_c)>-\frac 1 2 $ will be called ``end-pole'', and ``flag-pole'' if  {  $\Re(r_c)\leq -\frac{1}{2}$.  
}
\end{itemize}
\end{definition}

The motivation behind  this definition is that we can integrate the weight function  $w(z) = {\rm e}^{\theta(z)} \d z$ on contours that approach a  pole $c$ along the  steepest descent directions (if $d_c\geq 1$) or along any direction in the case of {end}-poles  ($d_c=0, \ \Re r_c>-\frac{1}{2}$) while obtaining a well defined integrable integrand.  Consequently we will denote by 
$
\int_{c^{(\ell)}}^{ s^{(m)} }
$
 an integration along a path that approaches the poles $c$ and $ s$ along the specified directions $c^{(\ell)}$ and $ s^{(m)} $.

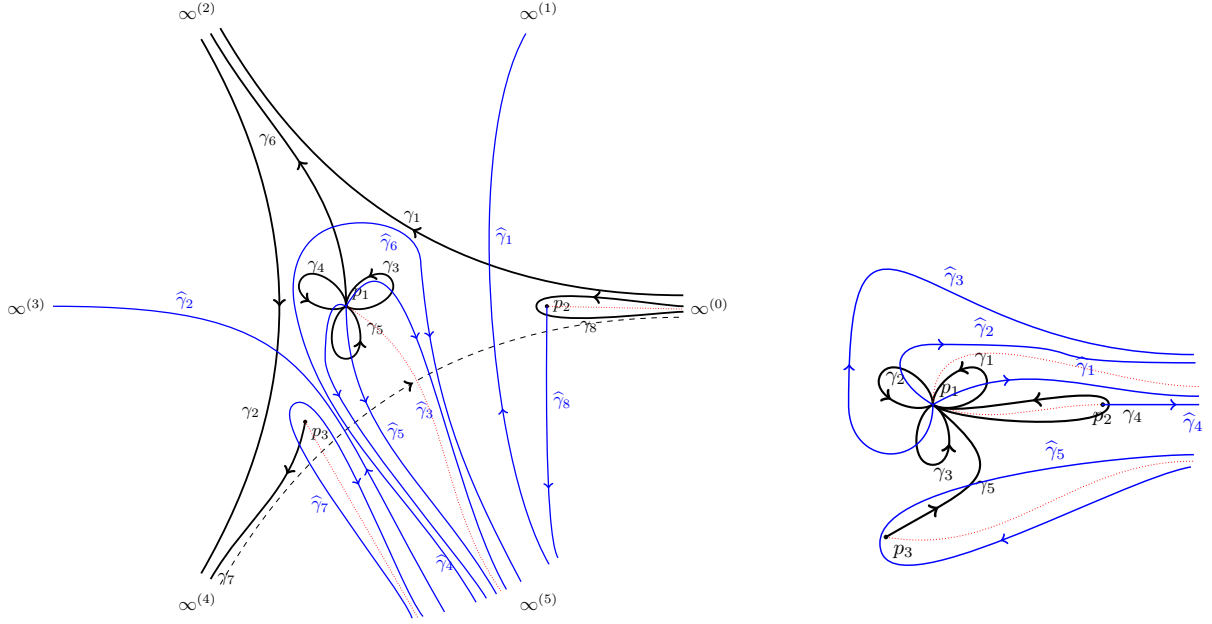
\begin{figure}[t]
\begin{center}
\resizebox{0.6\textwidth}{!}{
\begin{tikzpicture}
\draw [line width=1, postaction={decorate,decoration={markings,mark=at position 0.5 with {\arrow[line width=1.5pt]{>}}}}] (2:6) to  [out=180, in=-60]node[above] {$\gamma_1$} (118:6);
\draw [line width=1, postaction={decorate,decoration={markings,mark=at position 0.5 with {\arrow[line width=1.5pt]{>}}}}] (122:6) to [out=-60, in=60] node[pos=0.7,left] {$\gamma_2$}(238:6);
\draw [dashed, line width=0.2, postaction={decorate,decoration={markings,mark=at position 0.5 with {\arrow[line width=1.5pt]{>}}}}] (242:6) to [out=60, in=180] (-2:6);

\begin{scope}[xshift=-12]
\draw [line width=1, postaction={decorate,decoration={markings,mark=at position 0.7 with {\arrow[line width=1.5pt]{>}}}}] 
(0,0) to [out=-30, in=-60] node[pos=1,above] {$\gamma_3$} (30:1) to  [out=120,in=90]
 cycle;
 \draw [fill] (0,0) circle [radius=0.84pt]  node [above right] {$p_1$};
 
\draw [line width=1, postaction={decorate,decoration={markings,mark=at position 0.7 with {\arrow[line width=1.5pt]{>}}}}] 
(0,0) to [out=90, in=60] node[pos=1,above right] {$\gamma_4$} (150:1) to  [out=240,in=-150]
 cycle;
\draw [line width=1, postaction={decorate,decoration={markings,mark=at position 0.7 with {\arrow[line width=1.5pt]{>}}}}] 
(0,0) to [out=-150, in=180]  (-90:1) to  [out=0,in=-30] node[pos=0.6,right]{$\gamma_5$} 
 cycle;
 
\end{scope}

\draw [line width=1, postaction={decorate,decoration={markings,mark=at position 0.5 with {\arrow[line width=1.5pt]{>}}}}] (-0.42,0) to [out=90, in=-60] node[pos=0.6,left] {$\gamma_6$} (120:6) ;
\draw [red, densely dotted] (-0.42,0) to [out=-30, in=125]  (-65:6) ;
\draw [line width=1, postaction={decorate,decoration={markings,mark=at position 0.3 with {\arrow[line width=1.5pt]{>}}}}] (0:6) to[out=180, in =90, looseness=0.4] (3.2,0)  to[out=-90, in=-180, looseness=0.4] (-1:6) ;
\node at(4.2,-0.4) {$\gamma_8$};
\draw [fill] (3.4,0) circle[radius=0.84pt] node[right]{$p_2$};
\draw [red, densely dotted] (3.4,0) to (-0.5:6);
\draw [fill] (-1.2,-2.2) circle[radius=0.84pt] node[below right] {$p_3$};
\draw [line width=1, postaction={decorate,decoration={markings,mark=at position 0.3 with {\arrow[line width=1.5pt]{>}}}}] (-1.2,-2.2) to [out=-100,in=60] node[pos=1, right] {$\gamma_7$} (-120:6) ;
\draw [red, densely dotted] (-1.2,-2.2) to [out=-60,in=120]  (-81:6) ;
\draw [blue,line width=0.7, postaction={decorate,decoration={markings,mark=at position 0.3 with {\arrow[line width=1pt]{<}}}}] (-53:6) [out=120,in=-90, looseness=0.3] to node[pos=0.6, right] {$\wh \gamma_8$} (3.4,0);

\draw [blue,line width=0.7, postaction={decorate,decoration={markings,mark=at position 0.3 with {\arrow[line width=1pt]{>}}}}] (-55:6) [out=120,in=240, looseness=0.6] to node[pos=0.6, right] {$\wh \gamma_1$} (60:6);
\draw [blue,line width=0.7, postaction={decorate,decoration={markings,mark=at position 0.3 with {\arrow[line width=1pt]{>}}}}] (-76:6) [out=120,in=0, looseness=1.4] to node[pos=0.8,above] {$\wh \gamma_2$} (180:6);

\draw [blue,line width=0.7, postaction={decorate,decoration={markings,mark=at position 0.3 with {\arrow[line width=1pt]{>}}}}] (-0.42,0) [out=60,in=120, looseness=1] to node[pos=0.6,left] {$\wh \gamma_3$} (-64:6);

\draw [blue,line width=0.7, postaction={decorate,decoration={markings,mark=at position 0.3 with {\arrow[line width=1pt]{>}}}}] (-0.42,0) to [out=150,in=90, looseness=1] (-0.82,-1)
to [out=-90,in=120, looseness=0.4] 
 node[pos=0.8,left] {$\wh \gamma_4$} (-68:6);

\draw [blue,line width=0.7, postaction={decorate,decoration={markings,mark=at position 0.3 with {\arrow[line width=1pt]{>}}}}] (-0.42,0) to [out=-90,in=120, looseness=1] 
 node[pos=0.4,right] {$\wh \gamma_5$} (-66:6);

\draw [blue,line width=0.7, postaction={decorate,decoration={markings,mark=at position 0.3 with {\arrow[line width=1pt]{<}}}}]  (-80:6) to [out=130,in=30,looseness=0.4]  ($(-1.2,-2.2)+(120:0.4)$)
to[out=210, in=110, looseness=0.4] node[pos=0.5, left] {$\wh \gamma_7$} (-82:6);

\draw [blue,line width=0.7, postaction={decorate,decoration={markings,mark=at position 0.3 with {\arrow[line width=1pt]{<}}}}] (-61:6) to[out=118, in=-85] (1,0.8)
to[out=90,in=19]node[below]{$\wh \gamma _6$}(-0.7,1.5) to[out=200, in=100] (-1.4,0) to[out=-80, in=120, looseness=0.9] (-70:6);

\foreach\a in{0,1,2,3,4,5}
{
\node at (\a*60:6.5) {$\infty^{(\a)}$};
}
\end{tikzpicture}}\qquad
\qquad
\resizebox{0.3\textwidth}{!}{
\begin{tikzpicture}

\begin{scope}[xshift=-12]

\draw [line width=1, postaction={decorate,decoration={markings,mark=at position 0.7 with {\arrow[line width=1.5pt]{>}}}}] 
(0,0) to [out=-30, in=-60] node[pos=1,above] {$\gamma_1$} (30:1) to  [out=120,in=90]
 cycle;
  \draw [fill] (0,0) circle [radius=0.84pt]  node [above right] {$p_1$};

\draw [line width=1, postaction={decorate,decoration={markings,mark=at position 0.7 with {\arrow[line width=1.5pt]{>}}}}] 
(0,0) to [out=90, in=60] node[pos=1.2,above right] {$\gamma_2$} (150:1) to  [out=240,in=-150]
 cycle;
\draw [line width=1, postaction={decorate,decoration={markings,mark=at position 0.7 with {\arrow[line width=1.5pt]{>}}}}] 
(0,0) to [out=-150, in=180]  (-90:1) to  [out=0,in=-30] node[pos=0.2,below]{$\gamma_3$} 
 cycle; 
\end{scope}
\draw [red, densely dotted] (-0.42,0) to [out=-20, in =180](2.4,0)  ;
\draw [red, densely dotted] (-0.42,0) to [out=90, in =180](4,0.3)  ;
\draw [red, densely dotted] (-1.2,-2.2) to [out=-10, in =180](3.9,-0.94)  ;


\draw [line width=1, postaction={decorate,decoration={markings,mark=at position 0.3 with {\arrow[line width=1.5pt]{<}}}}] (-0.42,0) to[out=-40, in =90, looseness=0.4] (2.5,0)  to[out=-90, in=-41, looseness=0.4] (-0.42,0) ;

\node at(2.9,-0.2) {$\gamma_4$};
\draw [fill] (2.4,0) circle[radius=0.84pt] node[below ]{$p_2$};

\draw [fill] (-1.2,-2.2) circle[radius=0.84pt] node[below right]{$p_3$};
\draw [line width=1, postaction={decorate,decoration={markings,mark=at position 0.3 with {\arrow[line width=1.5pt]{>}}}}] (-1.2,-2.2) to [out=30,in=-40,looseness=2] node[pos=0.4, right] {$\gamma_5$} (-0.42,0) ;

\draw [blue,line width=0.7, postaction={decorate,decoration={markings,mark=at position 0.3 with {\arrow[line width=1pt]{<}}}}] (0:4) [out=180,in=180, looseness=0.3] to node[pos=0.1, below] {$\wh \gamma_4$} (2.4,0.0);

\draw [blue,line width=0.7, postaction={decorate,decoration={markings,mark=at position 0.3 with {\arrow[line width=1pt]{>}}}}] (-0.42,0) [out=30,in=180, looseness=1] to node[pos=0.6,above] {$\wh \gamma_1$} (2:4);

\draw [blue,line width=0.7, postaction={decorate,decoration={markings,mark=at position 0.3 with {\arrow[line width=1pt]{>}}}}] (-0.42,0) [out=150,in=180, looseness=2] to ++(0,1) to [out=0,in=180] node[pos=0.1,above] {$\wh \gamma_2$} (10:4);

\draw [blue,line width=0.7, postaction={decorate,decoration={markings,mark=at position 0.3 with {\arrow[line width=1pt]{>}}}}] (-0.42,0) [out=-90,in=-90, looseness=2] to ++(-1.4,0) to [out=90,in=180] node[pos=0.6,above] {$\wh \gamma_3$} (12:4);

\draw [blue,line width=0.7, postaction={decorate,decoration={markings,mark=at position 0.3 with {\arrow[line width=1pt]{>}}}}] (-15:4) to  [out=190,in=-90, looseness=0.6]  (-1.3,-2.3) to [out=90,in=180,looseness=0.5] node[pos=0.6,above] {$\wh \gamma_5$} (-12:4);


\end{tikzpicture}}
\end{center}
%
%
%
\caption{Left: $B(z)= (z-p_1)^4 (z-p_2)(z-p_3)$ and $A$ is a relatively prime polynomial of degree $8$  with positive leading coefficient  so that $d=8$.  In black the contours  $\gamma_1,\dots,\gamma_8$  for  the case  $d_\infty=3$ and $d_1 = 3$, (the center of the ``shamrock'', $p_1$), $d_2=d_3=0$. We also have a {flag}-pole  (the point surrounded by the ``lasso'', $p_2$ ) and an {end}-pole (the endpoint $p_3$ of $\gamma_7$ ). 
{\color{black}
The dual contours go to infinity along the direction $\omega^{(5)}$. In this example  {the dotted red arcs denote the branch-cuts of ${\rm e}^{\theta(z)}$ }.
}
Right: the case where $B$ is as above but $A$ is a polynomial of degree less than that of $B$.  {\color{black} The dual contours go to infinity in an arbitrary direction}. 
}
\label{contours}
\end{figure}

For each pole $c \in \mathbb{P}^1$ we now describe certain contours emanating from it. 
We give first the description of these contours  for the case $\deg A\geq \deg B$ so that  $d_\infty  \geq 1$ and  $U_\infty$ is a polynomial of degree $\geq 1$.
\begin{definition}[Contours $\gamma_j$] \label{contoursDEF}
For each pole $c \in \mathbb{P}^1$ of $\theta'(z) \d z$ we define a number of contours depending on the order of the pole.
\begin{itemize}
\item 
For each pole $c \neq \infty$ of order $d_c+1\geq 2$ we choose $d_c$ contours (``petals'') approaching $c$ along the  consecutive steepest descent directions   $c^{(2k)}$,   $c^{(2k+2)}$.  We also pick a contour (``stem'') extending from $c^{(0)}$ to a steepest descent direction at $\infty$. 
\item
For the pole at $c=\infty$ {  of order $d_\infty\geq 2$} we choose  $d_\infty-1$  contours between $[\infty^{(2k)}, \infty^{(2k+2)}]$, $k=0,\dots,d_\infty-2$, with the definition of the steepest descent directions as in Def. \ref{directions}. Note that there is no contour between  $\infty^{(2d_\infty-2)}$ and $\infty^{(0)}$
 (see Fig. \ref{contours})\footnote{\color{black} Note that for the Freud weight $A(z)=z^a$  we have $d_{\infty}=a+1=d+1$.}. These contours will be taken so as to leave all zeros of $B$ on their left region.
 {  In the case $d_\infty=1$  there is no contour.  }
\item
For each {end}-pole $c$ (including the hard-edges, i.e. the zeros of $B$ that are not poles of $\theta')$ we pick a contour from $c$  to $\infty$ along a steepest descent direction.

\item
For each {flag}-pole $c$ with non-integer residue we choose a contour (``lasso'') coming from $\infty^{(2k)}$, circling $c$ in the counter-clockwise direction and returning to $\infty^{(2k)}$, where $k$ is any choice.
If the {flag}-pole has negative integer residue, we can replace this choice by a small circle. 
\end{itemize}
It is understood that the contours are chosen by avoiding the zeros of $B$ except possibly at the endpoints,  and are chosen so that they do not intersect each others except, possibly, at endpoints.  {Moreover, the branch-cuts, denoted collectively by $\scr B$,  of the function ${\rm e}^{\theta(z)}$ are chosen as follows (see Fig. \ref{contours} for illustration):} 
\begin{enumerate}
\item[-] for each pole $c\neq \infty$ of order $d_c+1\geq 2$ the cut extends to infinity along the direction $\infty^{(2d_\infty-1)}$;
\item [-]  the same for end-poles ;
\item [-] for a flag-pole, the branch-cut is chosen to extend to infinity within the lasso.
\end{enumerate}

\end{definition}

\paragraph{Dual contours and intersection pairing.}
For reasons that will become apparent but that already motivated a similar construction in \cite{Bertosemiclassical}, we need to define ``dual contours'' and a notion of intersection pairing. 

To give a sense of the motivation, we mention  that the space of semiclassical moment functionals of type $(A,B)$ is in duality with that of type $(-A-B', B)$: the weight functions are 
\be
(A,B)&\mapsto {\rm e}^{\theta(z)}, \ \ \ \theta'(z) = -\frac {A(z)+B'(z)}{B(z)}\\
(-A-B',B)&\mapsto  \frac{{\rm e}^{-\theta(z)}}{B(z)} = {\rm e}^{\wh \theta(z)}.
\label{dualsym}
\ee
This duality maps the steepest descent directions of one symbol function into the ascent directions of the other, the {end}-poles (including hard-edges) of one into the flag-poles of the other, and viceversa.\footnote{Since we have stipulated that the  poles $c$ with $r_c=-\frac 1 2$ are considered {flag}-poles for the functional $\mathcal M$, then for the dual functional  they will be treated  as {end}-poles.}
The simplest explanation of the duality is by considering a generating function $\Phi(\l) =  \mathcal M[{\rm e}^{\l z}]$ of type $(A,B)$ and a generating function $\Psi(\l) =  \wh{\mathcal M}[{\rm e}^{-\l z}]$ of the dual type {\color{black} defined similarly to $\mathcal M$ but with the symbol  $\wh \theta$ \eqref{dualsym} and integration along suitable contours (Def. \ref{defdualcont}).} A simple exercise shows that they satisfy the two adjoint differential equations:
\be
\Big(\l B(\pa_\l) -A(\partial_\l) \Big)\Phi(\lambda)=0 \qquad
\Big( B(-\pa_\l)\l -A(-\pa_\l)\Big)\Psi=0.  \label{ODEgen}\ee
The solution spaces of two adjoint equations are put in duality by the {\it bilinear concomitant} \cite{Ince}: for equations with linear coefficients in $\l$ the bilinear concomitant has a homological interpretation as intersection pairing of the dual contours that we are describing here {\color{black} and are illustrated in Fig. \ref{contours}.}

\begin{definition}[Dual contours $\widehat{\gamma}_j$]
\label{defdualcont}
For each pole $c \in \mathbb{P}^1$ we define a number of contours dual to the contours $\gamma_j$ in Def.~\ref{contoursDEF}  {so as to avoid the branch-cuts $\scr B$:}
\begin{itemize}
\item
For each pole $c$ of $\wh \theta' \d z$ order $d_c+1\geq 2$ we choose $d_c$ contours (``anti-petals'')   {\color{black} originating from $c$ along the directions $c^{(2k+1)}, k=0,\dots, d_c-1$   and extending to $\infty^{(2d_\infty-1)}$. }
We also pick a contour (``anti-stem'') lassoing $c$ from $\infty^{(2d_\infty-1)}$ and intersecting only the corresponding stem.

\item
For the pole at $c=\infty$ we choose   $d_\infty-1$  contours   $[ \infty^{(2d_\infty-1)},\infty^{(2k+1)}]$, $k=0,\dots, d_\infty-2$. Note that they can be arranged so as to intersect only the corresponding steepest descent contours (see Fig. \ref{contours}). 
\item
For each {end}-pole $c$ of the original functional (which is a {flag}-pole for the dual functional)  we pick a lasso from $\infty^{(2d_\infty-1)}$.
If the residue of $c$ in $\theta'\d z$ is a positive integer {\color{black} (i.e.  a zero of the weight ${\rm e}^{\theta}$ (end-pole) and a pole of the dual weight ${\rm e}^{-\theta}/B$)}, we will choose a circle.
\item
For each {flag}-pole $c$ we choose a contour $[c,\infty^{(2d_\infty-1)}]$ this includes the poles, if any, where $r_c = -\frac 12 $ (which were considered as {flag}-poles). 
\end{itemize}
\end{definition}

The construction is such that for each contour $\gamma_j$ there is exactly one dual contour $\wh \gamma_j$  that intersects $\gamma_j$ at a single point $\gamma_j$ and has no intersections with any other contour. The orientations are chosen so that 
\be
\wh \gamma_j \circ  \gamma_k = \delta_{jk},\label{intpair}
\ee

%
\subsubsection{The case $\deg A<\deg B$ and $B$ with at least one multiple root} 
\label{BA}
In this case we 
 assume that $B$ has at least one zero, $z=p_0$,   of multiplicity higher than $1$. 
Note that in this case the condition that $A$ and $B$ are relatively prime  prevents $A$ to be a multiple of $B'$.
See Remark \ref{finiterank}  and Appendix~\ref{family}   {for a discussion of the case where $A$ is an integer multiple of $B'$.} 
 Then the symbol $\theta$ has  a pole at $z=p_0$ with at least one steepest descent direction and one steepest ascent direction. 
We keep the same terminology  { as in the previous case}. The contours will be chosen similarly as in the previous case but with the  ``stems'', ''lassoes'', and the contours to the flag-poles extending to $p_0^{(1)}$ (the point $p_0$ along the first-steepest ascent direction) instead of $\infty^{(2d_\infty-1)}$. 
The dual contours  are exactly as before, extending to $\infty$ (in any direction, for example the positive axis). The reason for this definition of dual contours in this case is motivated by the use of them  we need to make  in the main Theorems \ref{thm1}, \ref{thm2}.
\begin{remark}
\label{remark2.5}
We will not treat the case where all the zeros of $B$ are simple  (the classical Heine-Stieltjes electrostatic problem); in this case 
 the contours should be chosen,  generically (i.e. for non-integer residues of $\theta'$) as Pochhammer contours (group-commutators of the generators in the fundamental group of the plane minus the zeros of $B$) but there are  case distinctions according to whether the residues of $\theta'$ at these poles are positive integers, negative integers or neither which  complicate the description.
 \end{remark}

{\color{black}
\br
\label{finiterank}
The case when $\deg A< \deg B$ allows in general the possibility that $A$ is a multiple of $B'$. If $A$ is an {\it integer} multiple of $B'$, $A = k B'$ with $k=1,2,3\dots $ then  the weight of the semiclassical moment functional  is 
\be
A = k B' \ \ \Rightarrow \ \ {\rm e}^{\theta(z)} = \frac 1{ B(z)^{k+1}}.
\ee
In this case the moment functionals $\mathcal M$ are of finite rank (they are linear combinations of derivatives of Dirac delta functions supported at the zeroes of $B$, with the order of derivative being equal to $k$ times the multiplicity of said zero). 

It then follows that the orthogonal polynomials of degree $n\geq (k+1) \deg B$ are all the polynomials divisible by $B^{k+1}$, a type of solution that we could call ``improper''. For $n= k\deg B + 1$  we find interesting, albeit simple, {\it continuous families} of proper solution of the Stieltjes--Bethe equations, see Appendix \ref{family}.
\er
}

\paragraph{Connection with Painlev\'e\ equations.}
These types of moment functionals were analyzed in \cite{Bertola:Semiiso}   where it was shown that the Hankel determinant of the moments, when considered as function of the coefficients of $A, B$, are ``isomonodromic tau functions'' in the sense of \cite{JMU}. In particular this means that specializing the symbol $\theta$ one can  connect the theory of orthogonal polynomials to certain solutions of the   Painlev\'e\ equations II,$\dots$, VI, as well as many integrable generalizations thereof. 

In this perspective the maximal degeneration of a polynomial implies   that the tau-function must vanish and hence, in the cases that overlap with the theory of Painlev\'e\ transcendents, we are considering poles of the corresponding transcendent. 

This connection was exploited in our companion paper \cite{BertoGravaHeredia} in connection with a conjecture of Shapiro and Tater \cite{ST22}.   

\section{Degenerate orthogonal polynomials and Fekete problem}
A polynomial $P_n(z)$ of degree $\leq n$ is {\it orthogonal} {\color{black} for a moment functional $\mathcal M$ if \cite{Chihara}}
\be
\langle P_n, z^k\rangle := \mathcal M\Big[P_n(z) z^k\Big] = 0,\  \ \ k=0,\dots, n-1.\nn
\ee
Consider now  semiclassical moment functional $\mathcal M$ of type $(A,B)$ and let $d=\max\{a,b-1\}$: then the orthogonality reads 
\be
\label{ortho1}
\mathcal M[P_n (z) z^k] = 
\sum_{j=1}^d s_j \int_{\gamma_j} P_n(z)z^{k} {\rm e}^{\theta} \d z=\int_{\Gamma} P_n(z)z^{k} {\rm e}^{\theta} \d z=0 ,\ \ \ k=0,\dots, n-1,
\ee
\noindent where $\Gamma= \sum_{j=1}^{d} s_j \gamma_j$ and $\theta$ is determined from $A,B$ by \eqref{semidef} (up to an inessential additive constant). This is sometimes called ``non-Hermitian'' orthogonality because the inner product is complex and bilinear rather than sesquilinear. 
Let $\mu_k({\bf s}) = \int_{\Gamma} z^k {\rm e}^{\theta}\d z$ be the moments and define 
\be
\label{Dn}
D_n({\bf s}):= \det \le[\mu_{j+k}({\bf s})\ri]_{j,k=0}^{n-1}\,.
\ee
The polynomials of degree $\leq n$ that satisfy the orthogonality \eqref{ortho1} are expressible as the determinant 
\be
\label{Pn}
P_n(z) =\det \le[
\begin{array}{cccc}
\mu_0 & \mu_1 & \dots & \mu_{n}\\
\mu_1 & \mu_2 & \dots &\mu_{n+1}
\\
\vdots
\\
\mu_{n-1} & \dots & &\mu_{2n-1}\\
1 & z & \dots & z^n
\end{array}
\ri]\,.
\ee
 We now recall the Definition~\ref{defdeg}.
The polynomial $P_n$ is called {\bf $\ell$--degenerate orthogonal} if, in addition to the orthogonality conditions \eqref{ortho1}, it satisfies 
\be
\le<P_n,z^{k}\ri>=0, \ \ \ k=0,1,\dots, n+\ell-1.
\ee
As we mentioned in the introduction, the notion of $\ell$--degeneracy translates into the vanishing of suitable determinants, as  shown in the next Lemma, whose proof is immediate.
\bl
\label{lemmaell}
The orthogonal polynomial $P_n(z)$ is $\ell$--degenerate  ($\ell\geq 1  $)    if and only if the following determinants vanish:
\be
D_{n+1,k}({\bf s}):= \det H_{n+1,k}=0,\ \ \ \ k=0,1,\dots, \ell-1,
\label{Dnk}
\ee
where $H_{n+1,k}$ are the matrices
\be
 H_{n+1,k}:= \le[
\begin{array}{cccc}
\mu_0 & \mu_1 & \dots & \mu_{n}\\
\mu_1 & \mu_2 & \dots &\mu_{n+1}
\\
\vdots
\\
\mu_{n-1} & \dots & &\mu_{2n-1}\\
\hline
\mu_{n+k}& \mu_{n+k+1}  & \dots & \mu_{2n+k}
\end{array}
\ri].
\ee
\el
 According to Definition \ref{defdeg}, a $0$--degenerate polynomial is just an orthogonal polynomial (no conditions are imposed) and  generically it  exists.  
 We will say that $P_n$ is {\bf maximally} degenerate if it is $d-1$--degenerate. 
 The justification of the terminology is that the condition of $\ell$--degeneracy imposes $\ell$ homogeneous polynomial equation constraints on the parameters $s_1, \dots, s_d$ and hence, generically,  we can impose at most $d-1$ such constraints while expecting to have solutions. 

{\color{black}
It is important to point out the following.  Suppose that $D_{n,0}({\bf s})=0=D_{n+1,0}({\bf s})$: one then  concludes that on this locus the expression \eqref{Pn} for the orthogonal polynomial $P_n(z)$ yields the identically zero polynomial. This is true for an arbitrary moment functional (namely in the indeterminates $\{\mu_\ell\}_{\ell\in \N}$) as we now show, starting from the following proposition.

\bp
\label{caustic}
Suppose that $D_{n,0} = D_{n+1,0}=0$. Then $P_n(z) \equiv 0$; more specifically,  considering $P_n(z)\in \C[z] \otimes \C[\mu_0,\dots, \mu_{2n}] $ then $P_n(z)^2$ lies in the ideal generated by $D_{n,0}$ and $D_{n+1,0}$.   In particular  we also have  $D_{n+1,k}=0, \ \forall k\in \N$.
\ep
\noindent 
{\bf Proof.}
Consider the expression for  $P_n(z)$ :
 \be P_n(z)=
\det\le[
\begin{array}{cccc}
\mu_0 & \mu_1 & \dots \mu_{n-1}& \mu_{n}\\
\mu_1 & \mu_2 &  & \vdots\\
\vdots &&\\
\mu_{n-1} & \mu_n &\dots \mu_{2n-2} &\mu_{2n-1}\\
1& z & \dots & z^n
\end{array}
\ri]=C_{n} z^n - C_{n-1} z^{n-1} + \dots + (-1)^n C_0
\ee
We are going to show that $C_j C_rk \in \le \langle D_{n,0}, D_{n+1,0}
\ri \rangle$ (the ideal generated by the two polynomials) for all $j,k=0,\dots, n$; this immediately implies the two subsequent statements.

Recall the (general) Desnanot-Jacobi identity; if $M$ is an $n\times n$ matrix and $J,K$ are two subsets of $\{1,\dots, n\}$ of the same cardinality $r$ (and listed in increasing order)  then we denote by $M^{J;K}$ the  $(n-r)\times(n- r)$ matrix obtained by deleting the rows indexed by $J$ and the columns indexed by $K$. 
Then the identity reads
\be
\label{desnanot1}
\det M \det M^{\{a,b\}; \{c,d\}} = \det M^{\{a\};\{c\}} \det M^{\{b\};\{d\}} - \det M^{\{b\};\{c\}} \det M^{\{a\};\{d\}},
\ee
{ for $a,b,c,d\in\{1,\dots,n\}$ with $a\neq b$ and $c\neq d$.}
Let $H$ be the  Hankel matrix of the moments of size $n+1$; then from \eqref{desnanot1} we obtain
\bea
\label{desnanot2}
\det H \det H^{\{j,n+1\}; \{k,n+1\}} =\det H^{\{j\};\{k\}} \det H^{\{n+1\};\{n+1\}} - \det H^{\{j\};\{n+1\}} \det H^{\{n+1\};\{k\}}.
\eea
The last term in the right side of \eqref{desnanot2} is precisely $C_{j}C_{k} $  {up to an overall  sign,}
so that we have 
\be
C_{j} C_{k} = \le(\det H^{\{j\};\{k\}}\ri) D_{n,0} -\le(\det H^{\{j,n+1\}; \{k,n+1\}}\ri) D_{n+1,0}  
\ee
This proves that all the quadratic expressions belong to the claimed ideal. \QED
}

 The relevance of the notion of $\ell$--degeneracy is clarified by the following proposition. 
\bp
\label{Wronsk}
Let $\mathcal M$ be a semiclassical moment functional of type $(A,B)$  { with the polynomials $A$ and $B$ relatively prime and
with $B$ that has at least one multiple root  when $\deg A<\deg B$.} Let $\theta$ be its symbol according to \eqref{semidef}. 
Given any polynomial $P_n(z)$ we set
\be
F(z):= \sqrt{B(z)}P_n(z){\rm e}^{\frac 1 2\theta(z)},\ \ \  G(z):= \sqrt{B(z)} R_n(z){\rm e}^{-\frac 1 2 \theta(z)}
\label{FG}\\
R_n(z):=  \frac{1}{2i\pi}\int_{\Gamma} \frac {P_n(x) {\rm e}^{\theta(x)}\d x}{(x-z)},
\ee
where  $\Gamma = \sum_{j}^d s_j \gamma_j$  and the contours $\gamma_j $ have been defined in Section~\ref{contourssec}.

The Wronskian $W= W(F,G) = F'G-G'F$ is a polynomial of degree $\leq \max\{\deg A+n-1, \deg B+n-2\}$.  If $P_n$ is an $\ell$--degenerate orthogonal polynomial, then $W$ is a polynomial of degree $d-\ell-1$, with $d=\max\{\deg A,\deg B-1\}$. 
{   If $P_n$ is maximally degenerate, then the Wronskian $W$ is equal to a constant.}
\ep
\noindent {\bf Proof.} We first show that $W$ does not have jump discontinuities across the contours $\gamma_j$ and extends to an entire function.
A direct computation shows:
\be
W = \theta' B P_n R_n + B(P_n' R_n-P_nR_n') =- \wh A  P_n R_n + B(P_n' R_n-P_nR_n') ,\label{Wdef}
\ee
where, for brevity, we have set $\wh A(z) = A(z)+B'(z)$.
Now, let $\Gamma = \sum_{j}^d s_j \gamma_j$ and take $z\in \gamma_j$. The remainder $R_n$ satisfies
\be
R_n(z_+) - R_n(z_-) &= P_n(z) {\rm e}^{\theta(z)},\\
R_n'(z_+) - R_n'(z_-)& = \Big(P_n'(z)  + \theta'(z) P_n(z)\Big) {\rm e}^{\theta(z)}
\ee
where we approach $z$ in the oriented contour $\gamma_j$ from the $\pm$ sides correspondingly.
Thus, with $\Delta R_n = R_n(z_+)-R_n(z_-)$ the jump operator, we have 
\be
W(z_+)- W(z_-) =& -\wh A  P_n \Delta R _n + B(P_n '\Delta R_n-P_n\Delta R_n')  
=\\
=&( BP_n' -\wh A  P_n ) \Delta R _n - B P_n \Delta R_n' 
=\\
=& ( BP_n'-\wh A  P_n ) P_n {\rm e}^{\theta }  - B P_n \Big(P_n'  + \theta' P_n\Big) {\rm e}^{\theta}
=0,
\ee
where in the last simplification we have used the equation \eqref{semidef} for $\theta'$. 
This concludes the proof of the absence of discontinuities. 
Now consider a zero of $B$; these are the only possible singularities of $\theta'$ (but they may also be regular points of $\theta'$, see Rem. \ref{catch22}).
From \eqref{Wdef} it is clear that the only possible singularities are at the endpoints of the contours $\gamma_j$. If this is the case of a petal or stem, then the integrand tends to zero exponentially and hence the singularity of $R_n$ is at worst logarithmic.
If this is an {end}-pole then the integrand in the definition of $R_n$ behaves as $(w-c)^{r_c}$ and hence the Cauchy transform  has at worst growth bounded by $\max\{|z-c|^{\Re (r_c)} , |\ln |z-c||\}$ \cite{Gakhov}. This shows that {\it a priori}, the expression \eqref{Wdef} for the Wronskian may have at worst an isolated singularity at the zero $z=c$ with growth bounded by $|z-c|^{\Re (r_c)}$ (if $\Re (r_c) \in (-1,0)$). So it actually must have a removable singularity and $W$ extends analytically also at the zeroes of $B$.

Suppose now that $P_n$ is an $\ell$--degenerate polynomial. This means that the remainder term $R_n$ is of order $\mathcal O(z^{-n-1-\ell})$ as $|z|\to\infty$. Then, from \eqref{Wdef} we see that  we have  
\be
W(z) = \mathcal O(z^{\max\{\deg A-\ell-1, \deg B-2 - \ell\}})= \mathcal O(z^{d-\ell-1}),\quad \mbox{as $|z|\to\infty$}.
\ee
{ Form the above relation we deduce that if $P_n$ is maximally degenerate,   namely $\ell=d-1$, then the Wronkstian is a constant.}
This completes the proof.
\QED

In view of Prop. \ref{Wronsk} the Wronskian $W$ of $F,G$ in \eqref{FG} is a {\it constant} when 
  $P_n$ is a maximally degenerate orthogonal polynomial.
   This is the crucial property needed to prove half of the Theorem \ref{thmintro}.
We first need the following simple Lemma
\begin{lemma} \label{nopoles}
Consider the second order ODE 
\begin{equation}
  \frac {\d^2 y}{\d z^2} -V(z)y =0.
\end{equation}
Suppose that  $z=z_*$ is (possibly) a   singularity of the potential $V(z)$ and it is a {pole of order
at most $2$}. Assume that  it is an  {apparent} singularity, namely, there are two linearly independent  solutions to the ODE that are analytic at $z=z_*$.
Then, $V(z)$ is locally analytic at $z_*$.
\end{lemma}
The proof is a direct local analysis of the indicial equation and can be found in detail in Lemma 3.6 of \cite{BertoGravaHeredia}.
{\color{black} Next we prove the following  theorem which is the first  part of Theorem \ref{thmintro}.}

 \bt
 \label{thm1}
 {  Let $\mathcal M$ be  a  semiclassical  moment functional of type $(A,B)$ with the polynomials $A$ and $B$ relatively prime and
with $B$ that has at least one multiple root  when $\deg A<\deg B$.} Suppose that $\mathcal M$  admits a maximally degenerate polynomial $P_n$ of degree $n$.
   Then 
 \begin{enumerate}
 \item [{\bf (1)}] The function $F(z) = \sqrt{B(z)} P_n(z) {\rm e}^{\frac{1}{2}\theta(z)}{  = P_n(z){\rm e}^{-\tfrac 12 \wh \theta(z)}}$ solves the differential equation 
 \be
 \label{SL}
 F''(z) - V(z) F(z) =0\,,
 \ee
 where the potential $V(z)$ is a rational function with poles only at the zeros of $B(z)$ at most of twice the order of the form:
 \be
 V(z) =  
 -\frac 1 2\wh\theta'' + \frac 1 4(\wh\theta')^2 
 + \frac {Q}{B} , \ \ \ \deg Q \leq d-1.
 \ee
  If $B\equiv 1$ the potential $V$ is a polynomial of degree $2\deg A$.
 Equivalently the polynomial $P_n$ satisfies the  generalized Heine-Stieltjes  equation:
 \be
 \label{ODEPn}
  B(z) P_n'' - A(z)  P_n' - Q(z) P_n(z) = 0 .
 \ee

 \item [{\bf (2)}] Let  $(z_1, \dots, z_m)$ denote the roots of $P_n$  { that are   distinct  from those of $B$};  they  satisfy the system of equations
 \be
 \sum_{\ell\neq j\atop 1\leq \ell\leq m }\frac 1{ z_j-z_\ell}   = \frac {A(z_j) }{2B(z_j)}
 - \sum_{c: B(c)=0} \frac {\sharp_c}{z_j - c},\ \ \ j= 1,\dots, m
 \label{Fekete}
 \ee
 { where  $c$ is a common root of $P_n(z)$ and $B(z)$ and  the positive integer $\sharp_c$ denotes the multiplicity of the root $c$ of $B(z)$.}
  {
 Under the genericity assumption that $A(c) - kB'(c)\neq 0$ for all roots  $c$ of $B$ and all $k =0,1, \dots$, the polynomial $P_n$ does not share any root with $B$. }
 \end{enumerate}
 \et
\noindent{\bf Proof.}
{\bf (1)} By Proposition \ref{Wronsk} the Wronskian of $F, G$ defined in \eqref{FG} is a constant (and necessarily non-zero). 
Then we can write  (up to a rescaling)
\be
W = F G'- F' G=1 \ \Rightarrow \ W' = F G''- F'' G=0.\label{0W}
\ee
This means that $\frac {F''}F =  \frac {G''}G$ and we can recast the equation \eqref{0W} as a differential  equation of the form 
\be
\label{sturmunddrang}
y''(z) - V(z)y(z)=0,\ \ \ \ V := \frac {F''}F = \frac {G''}G.
\ee
The potential $V$ in \eqref{sturmunddrang} is {\it a priori}  a rational function with poles at all the zeros of $F$ (i.e. of $P$) as well as  singularities at the zeros of $B$.  
The essential part is to show that $V$ is analytic at  { those  zeros of $P$ that do not coincide with a zero of $B$ as well}. 

 To see this we observe that \eqref{sturmunddrang} has both $F$ and $G$ as solutions: the function $G$ as presented has discontinuities proportional to $F$ across the contours $\gamma_j$'s:    so we can analytically continue $G$ to the universal cover of the plane minus the zeros of $B$. 
Let $P_n(z) = \prod_{j=1}^n (z-z_j)$ {  and let us order them so that $z_1,\dots, z_m$, $m\leq n$ are the zeros distinct from those of $B$}; then both $F$ and (the above mentioned analytic continuation of) $G$ are   {analytic at} $z_j$, $j=1,\dots, {m}$,  which is therefore an {\it apparent singularity}\footnote{In the literature of Sturm--Liouville equations like \eqref{sturmunddrang} it is customary to call ``apparent singularity'' a pole of $V$ such that both solutions have Puiseux series expansion in half--integer powers. Here we use the terminology ``apparent singularity'' in the strict sense that both solutions of the equation must be locally analytic. } of the  equation \eqref{sturmunddrang}.  Then Lemma \ref{nopoles} says that $V$ must be locally analytic at $z_j$.
By definition of $V$ we have
\be
\label{Vexp}
V(z) =  \frac {F''(z)}{F(z)} =
 \frac 1 2\theta'' + \frac 1 4(\theta')^2 +\frac {B''}{2B} -\le(\frac {B'}{2B}\ri)^2+ \frac {B'}{2B} \theta' 
 + \frac {P_n''} {P_n} + \frac {\theta' P_n'}{P_n} + \frac {B'P_n'}{B P_n}.
\ee
By the Lemma \ref{nopoles} it follows that $V(z)$ given by \eqref{Vexp}  must be analytic at all  roots $z_j$, $j=1,\dots,  {m}$.

Since $V$ in \eqref{Vexp} has no poles at  {$z_1,\dots, z_m$},  the last three terms in \eqref{Vexp} are of the form 
\be
\frac {P_n''} {P_n} + \frac {\theta' P_n'}{P_n} + \frac {B'P_n'}{B P_n} = \frac {Q(z)}{B(z)}, 
\ee
with $Q(z)$ a {\it polynomial} of degree at most $d-1$. 
The ODE \eqref{ODEPn} for $P_n$ follows then by straightforward manipulations. 
\\[5pt] 
\noindent{\bf (2)}
Since we have established that $V$ is analytic at each $z_j$, {  $j=1,\dots, m$}, it then follows that these  must be simple zeros of $P_n$ since $F$ is a non-trivial solution to a second--order ODE. Then we express the fact that   $\res{z=z_j} V(z)\d z=0$  with $V$ given by \eqref{Vexp} and obtain 
\begin{equation}
\label{SF}
2\sum_{\ell\neq j\atop 1\leq \ell\leq n} \frac 1{z_j-z_\ell} = \res{z=z_j} \frac {P''_n(z)}{P_n(z)}  \d z = - \theta'(z_j) - \frac {B'(z_j)}{B(z_j)} = \frac {A(z_j)}{B (z_j)},\ \ \  {j=1,\dots, m}.
\end{equation}
{ 
Consider now those zeros of $P_n$, $z_{m+1},\dots, z_n$ that are shared with those of $B$ (in principle they could be repeated). 
We move the terms containing them  from the left side of \eqref{SF} to its right side:
\be
2\sum_{\ell\neq k\atop 1\leq\ell\leq m} \frac 1{z_j-z_\ell}= \frac {A(z_j)}{B (z_j)} - \sum_{\ell=m+1}^n \frac 1{z_j-z_\ell},\ \ \  {j=1,\dots, m}.
\ee
This proves formula \eqref{Fekete}. Finally, 
we prove the genericity statement that if $k\,B'(c) - A(c) \neq 0$ for all zeros $c$ of $B$ and for all $k=0,1,2,\dots$, then $P_n$ does not share any zeros with $B$.   {  Suppose by contradiction that $P_n$ and $B$ share a zero $c$. Then $P_n\equiv 0$.}    This is seen using   the fact that  $P_n$ solves the ODE \eqref{ODEPn}.
Indeed, differentiating the ODE  \eqref{ODEPn}  $k$ times we get
\begin{align}
\sum_{s=0}^k& \le({k\atop s}\ri)\le[B^{(k-s)} P_n^{(s+2)} -A^{(k-s)}P_n^{(s+1)} - Q^{(k-s)}P_n^{(s)}\ri]=
\nn\\
=&
B(z) P_n^{(k+2)}(z) +  \big(k\,B'(z) - A(z)\big)P_n^{(k+1)}(z) + \big\{\text{ linear combination of $P_n, P_n',\dots, P_n^{(k)}$}\big\} =0.
\end{align}
Evaluating this latter equation at $z=c$ yields a recurrence relation for $P_n^{(k+1)}(c)$ in terms of a linear combination of lower derivative terms (note that the term $P_n^{(k+2)}$ does not enter into this recurrence because it is multiplied by $B(c)=0$).
It then follows by induction on $k$ that $P_n^{(k+1)}(c) =0$ for all $k$ and hence $P_n$ would be identically zero. 
}
\QED
{\color{black}
\br
\label{Lamerem}
The equation \eqref{ODEPn} for $P_n$
$$
 B(z) P_n'' - A(z)  P_n' - Q(z) P_n(z) = 0 
$$
 is a (degenerate) Lam\'e\ equation in the terminology of \cite{Boris18} {  or alternatively the  generalized Heine-Stiltjes  equation. The classical  Heine-Stiltjes  equation  occurs when $\deg A<\deg B$ and $B$ has  only  simple roots.}
\er
}
We now prove the converse of the Theorem  \ref{thm1} (and the second  part of Theorem \ref{thmintro}) in the following form.
\bt
\label{thm2}
Suppose that, for two relatively prime polynomials $A(z)$ and $ B(z)$,  a  solution of the Stieltjes--Fekete equilibrium problem \eqref{Fekete} consists of $n$ points 
$z_1,\dots, z_n$ (necessarily distinct).
 In the case $\deg A < \deg B$ we make the assumption that the polynomial $B$ has at least one zero of higher multiplicity and that 
 \be
 \label{bound}
 2n > \Re \Lambda -1 - \deg B, \ \ \ \ \Lambda:= - \res{z=\infty} \theta'(z)\d z.
 \ee
Then the polynomial $P_n(z) = \prod_{k=1}^n (z-z_k)$ is a maximally degenerate orthogonal polynomial for the pairing \eqref{ortho1}, with  the parameters $s_j$ given by  
\be
\label{sj2}
s_j = \int_{\wh \gamma_j} \frac {{\rm e}^{-\theta(z)}}{B(z)P_n^2(z)} \frac { \d z}{2i\pi}.
\ee
Here $\wh\gamma_j$ is the dual path to $\gamma_j$ in the homology as defined in Section \ref{contourssec}.

\et
\noindent{\bf Proof.}
The proof is mostly a back-tracking of the proof of Theorem \ref{thm1}. 
First of all the condition \eqref{Fekete} is stating that the expression for $V(z)$ in \eqref{Vexp}  with $F(z) = \sqrt{B(z)} P_n(z) {\rm e}^{\frac{1}{2}\theta(z)}$,  yields an analytic expression at all zeros of $P_n$. 
Now, with $V$ given by \eqref{Vexp}   we are seeking the linearly independent solution of the differential equation:
\be
\label{ODE2}
y(z)''- V(z) y(z)=0,\ \ \ \ \ \ V(z) = \frac {F(z)''}{F(z)}.
\ee
Using Abel's theorem (stating that the Wronskian of two solutions of \eqref{ODE2} is a constant) we can write a  second linearly independent 
solution by a quadrature 
\be
\label{intG}
G_q(z) = F(z) \int_{q}^z \frac {\d w}{F(w)^2}.
\ee
The basepoint $q$ of integration in \eqref{intG} can be chosen arbitrarily and different choices of basepoint amount to adding to $G_q$ a multiple of $F$.  { In the expressions the  branch-cuts $\scr B$ of the function $F$ are as specified in Definition \ref{contoursDEF} and Fig. \ref{contours}.}

The differential $\frac{\d w}{F(w)^2}$ has double poles at the zeros of $P_n$ but no residues; this is precisely guaranteed by the Fekete equilibrium equations \eqref{Fekete}. Therefore the antiderivative has simple poles without logarithmic singularity. 
Upon multiplication by $F(z)$ (which has simple zeros) the poles will cancel and  this provides the proof that $G(z)$ is locally analytic at all zeros of $P_n$. 

{\bf Case $\deg A\geq \deg B$}. 
Consider the  connected components of $\C\setminus \Gamma = \sqcup_\mu \scr D_\mu$: in the region $\scr D_{0}$ that contains $\infty^{(2d_\infty-1)}$ we use the latter as basepoint of integration. According to our choice of contours $\gamma_j$ (Sec. \ref{contourssec}) in every other connected component $\scr D_\mu$, {\color{black}  such that the boundary is $\gamma_j= \pa \scr D_\mu$},  there is exactly one endpoint, denoted by $c^{(\mu)}$  of a dual contour $\wh \gamma_j$ {\color{black} (this implies possibly also a particular direction of approach)}; in those regions we therefore define $G_{c^{(\mu)}}(z)$ by using the {\color{black}  basepoint $c^{(\mu)}$ for } integration in the formula \eqref{intG}.
 Then  patching together  these functions $G_{c^{(\mu)}}$ in each component $\scr D_\mu$  we {\it define} a piecewise analytic function $G(z)$ by 
\be
G(z) = F(z) \int_{c^{(\mu)}}^z \frac {\d w}{F^2(w)}{\color{black} := G_{c^{(\mu)}}(z)}, \ \ \ z\in \scr D_\mu.
\ee
 {Observe that the regions inside a lasso (e.g. the one containing $p_2$ in Fig. \ref{contours}) are not simply connected: if $c$ denotes the end-pole (thus with $\Re(r_c)>-\frac 12$) inside a lasso, then $G_c(z)$ has a branch-cut (oriented towards infinity, see Fig. \ref{contours}) across which $G_c(z_+) ={\rm e}^{-i\pi r_c}G_c(z_-)$ and a local behaviour $G(z) = (z-c)^{-\frac {r_c}2+1}\mathcal O(1)$.}
If $\gamma_j$ is the boundary of $\scr D_\mu$,  elementary calculus shows that for $z\in \gamma_j$ we have 
\begin{equation}
\label{324}
G(z_+) -G(z_-) = F(z) \int_{c^{(\mu)}}^{\infty^{(2d_\infty-1)}}\hspace{-20pt}\frac {\d w}{F^2(w)}=
 F(z) \int_{\wh \gamma_j}\frac {\d w}{F^2(w)} = s_j F(z). 
\end{equation}
Now  define   
\be
\label{intR}
R_n(z) := \frac{G(z) {\rm e}^{\frac 1 2 \theta(z)}}{\sqrt {B(z)}} =P_n(z){\rm e}^{\theta(z)} \int_{c^{(\mu)}}^z \frac {{\rm e}^{-\theta(w)} \d w}{B(w)P_n^2(w)}, \ \ \ z\in \scr D_\mu,
\ee
 { and observe that in each region $\scr D_\mu$ the function  $R_n(z)$ has no branch-cuts along $\scr B$ because the branch-cuts of ${\rm e}^{\theta} $ eliminate the branch-cuts of the integral within $\scr D_\mu$. As a consequence of \eqref{324} }
 we have 
\be
\label{RHPR}
R_n(z_+) = R_n(z_-) + s_j P_n(z) {\rm e}^{\theta(z)}, \ \ \ z\in \gamma_j.
\ee
We can express $R_n$ alternatively as a Cauchy transform; indeed, the equation \eqref{RHPR} implies that 
\be
\label{Rcacuhy}
R_n(z) = H(z) + \sum_{j=1}^d s_j \int_{\gamma_j} \frac {P_n(w) {\rm e}^{\theta(w)}\d w}{2i\pi(w-z)},
\ee
for some {\it entire} function $H(z)$, which we now show to be identically zero. To see this consider the component $\scr D_k$ that contains the steepest ascent direction $\infty^{(2k+1)}$ then the integral representation \eqref{intR} shows that $R_n(z)$ is bounded by $\mathcal O(z^{-n-b+1}) $ within  {the sector extending from the direction $\infty^{(2k)}$ to $\infty^{(2k+2)}$ and also slightly further (by adjusting the directions of approach to $\infty$ of the boundary of $\scr D_k$)}.  Since this holds for all the other regions that contain the steepest ascent directions, we see that $R_n(z)$ decays like $z^{-n-b+1}$ in  { overlapping open sectors whose union contains } every direction. This implies that $H$ must be identically zero by  Liouville theorem, since the Cauchy integral is already bounded by $\mathcal O(z^{-1})$ a priori.

Now consider the Wronskian $W(G,F)=1$; using the decay of $R_n(z) = \mathcal O (z^{-1-\ell}), \ \ \ell\geq 0$ we deduce from  \eqref{Wdef} that  
\be
W(G,F) =- (A+B')  P_n R_n + B(P_n' R_n-P_nR_n')=  \mathcal O(z^{d+n-1 - \ell}) \quad \mbox{as $|z|\to\infty$}.
\ee
 Since $W\equiv 1$ we conclude that actually  $\ell = d+n-1$ and  $R_n(z) = \mathcal O(z^{-n-d})$. This implies that $P_n$ is maximally degenerate from inspection of the Cauchy integral representation \eqref{Rcacuhy} (with $H\equiv 0$).

\paragraph{The case $\deg A<\deg B$. }

 The only point in the proof where we have used the condition $\deg A\geq \deg B$ is in  the choice of the base-point of integration  $\infty^{(2d_\infty-1)}$,  as described in Section \ref{contourssec} for this case.
We now suppose that $\deg A<\deg B$ and that $B$ has at least one multiple zero, see Section \ref{BA}. 
The modifications to handle this case are  of minor nature; since  {now} none of the contours $\gamma_j$ extends to infinity, we can choose dual contours $\wh \gamma_k$ that have intersection $\delta _{jk}$ with the contours $\gamma_j$ and extend to infinity (in any direction) as explained in Section \ref{BA}.  Note that the only use we make of the dual contours is in the reasoning from  formula \eqref{324} to the end of the above proof. 

Now, the integrand in the expression  \eqref {intR} (or equivalent \eqref{sj2}) is integrable at infinity under the following conditions:
\be
-2n  + \Re \Lambda- \deg B  < -1 \ \ \Rightarrow \ \ \ 
2n > \Re \Lambda -1 - \deg B, \ \ \ \ \Lambda:= - \res{z=\infty} \theta'(z)\d z.
\ee
The condition is not automatically guaranteed for  small $n\in \N$ because $\Re \Lambda$ can be arbitrarily large. 
If this inequality holds, however, then the integrand is $\mathcal O(z^{\Lambda-2n-\deg B })$ and hence the integral with basepoint at infinity is  $\mathcal O(z^{\Lambda-2n-\deg B+1})$:  it then follows that the whole expression \eqref{intR} is $\mathcal O(z^{-n-\deg B +1 }) =\mathcal O(z^{-n-\deg B +1 }) = \mathcal O(z^{-n-d})$.  This implies the maximal degeneracy of $P_n$ from the expression \eqref{Rcacuhy}, where $H\equiv 0$ is established by the use of Liouville's theorem. 
We then  see that the rest of the reasoning as well as the proof of the expressions \eqref{sj2} proceed unimpeded.  
\QED

For the case $\deg A<\deg B$ in   Theorem  \ref{thm2} we have assumed  that $B$ has a double zero, as well as  the  condition \eqref{bound} on the number of points $n$. 
An interesting example with simple zeros of $B$ is reported in Appendix \ref{family}; in general it would not be very complicated to lift the condition on the multiple zero of $B$, the  only price being additional description of contours (Pochhammer  contours in the generic case).  However it is less clear how to lift  the bound  \eqref{bound} on the degree $n$ (i.e. the minimal number of points for our proofs to proceed unimpeded): we thus do not know whether the theorem fails in these circumstances or only the proof  needs  modification.

\br
\label{remBA2}

\begin{figure}
\begin{center}
{$
\begin{tikzpicture}[scale=1]
\draw [->, dashed] (-1,0) to (2,0);
\draw [->, dashed] (0,-1.5) to (0,1.5);
\draw  [line width=1, postaction={decorate,decoration={markings,mark=at position 0.4 with {\arrow[line width=1.5pt]{<}}}} ](0,0) to[out=180, in=180 ,looseness=2] (0,1) to [out=0,in=90,looseness=1.4] (1,0) to [out=-90,in=0,looseness=1.4] (0,-1) to [out=180, in=180,looseness=2] (0,0);
\draw [dashed] (0,0) to (-2,0);
\draw  [blue,line width=0.7, postaction={decorate,decoration={markings,mark=at position 0.7 with {\arrow[line width=1.5pt]{>}}}}] (0,0) to [out=0,in=180] (10:2);
\draw  [dashed, blue,line width=0.7] (10:2) to [out=0,in=180] ++(1,0);
\end{tikzpicture}
\atop 
\rule{0pt}{70pt}
$} \qquad \includegraphics[width=0.3\textwidth]{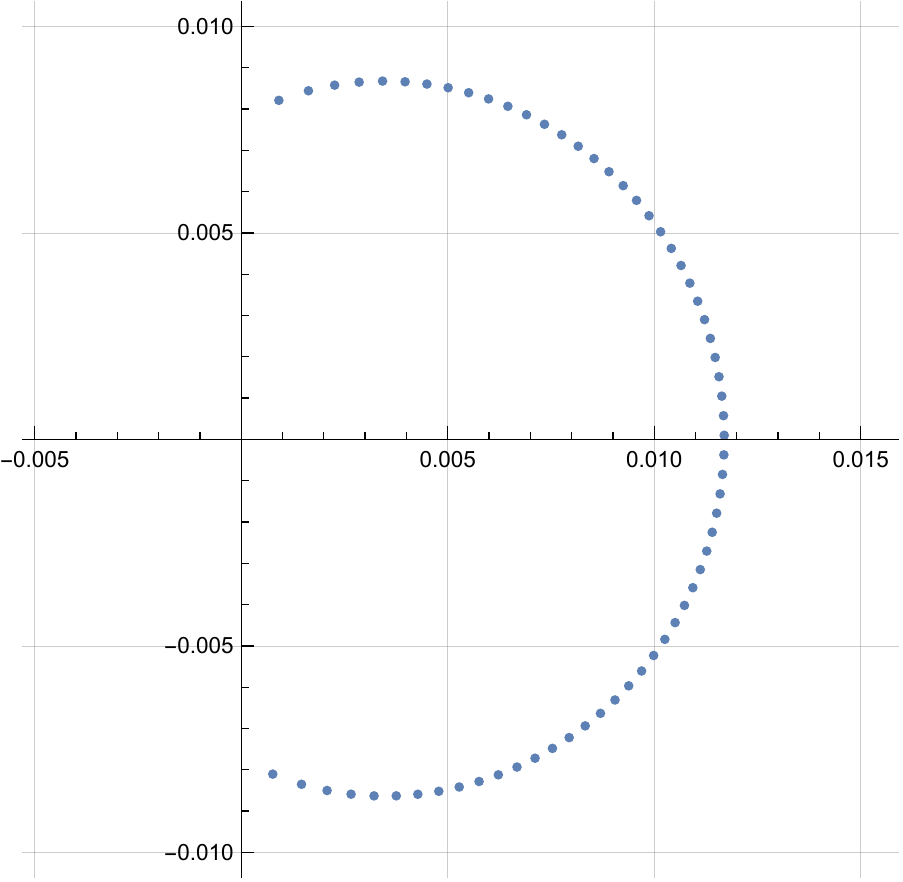}
\end{center}
\caption{Left: the contour and dual contour for the Bessel functional of Table \ref{classical}, for the case $\nu \not\in\Z$; for the case $\nu \in \Z$ one can choose simply a circle. Right: plot of $n=64$ points in critical configuration for the Bessel case with $A = 1-(2+i)z$, $B=z^2$. }
\label{figBessel}
\end{figure}
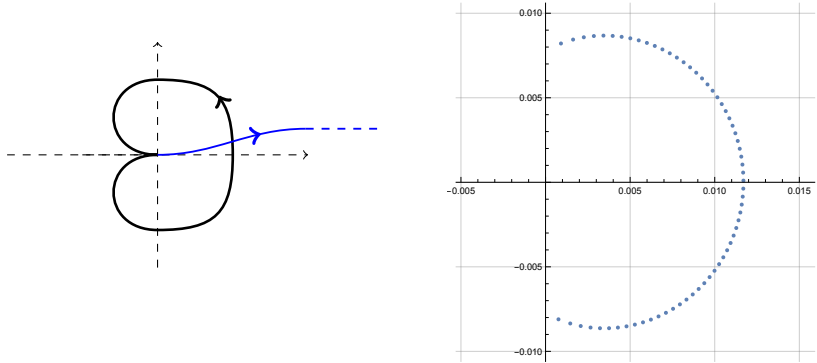

By the way of illustration let us consider the Bessel case in Table \ref{classical}; the only contour of integration is in Fig. \ref{figBessel} and as for dual we can take the positive real axis: the piecewise analytic remainder function $G(z)$ defined by  \eqref{FG}, outside of the cardioid coincides with the definite integral  $G_\infty(z):= F(z) \int_\infty^z \frac {\d x}{F^2(x)}$ (formula \eqref{intG}): here the direction of integration at $\infty$ is immaterial since $F(z)^2 \simeq z^{2n+\nu +2}$ (as long as $2n+\nu>-1$). Inside the cardioid $G(z)$ must be  recessive near $z=0$ along the positive direction and hence it must coincide with the definite integral  $G_{0^{(1)}} (z):= F(z) \int_{0^{(1)}} ^z \frac {\d x}{F^2(x)}$ , where we recall that the notation $0^{(1)}$ denotes the steepest ascent direction of the symbol $\theta$ near the point $z=0$. Thus we easily conclude that, for $z\in \gamma$ (with $\gamma$ the cardioid, oriented as indicated in Fig. \ref{figBessel}), $G(z_+) = G(z_-)  + s F(z) $ with $s = \int_0^\infty \frac {\d x}{F(x)^2}$, as claimed.
\er

\paragraph{Counting the number of solutions.}

A na\"ive counting would suggest, based on B\'ezout's theorem, that there are $(n+1)^{d-1}$ solutions of the set of equations that characterize the maximal degeneracy:
\be
\label{maxdeg}
D_{n+1, 0}({\bf s}) =0, \ \dots, \ \ D_{n+1,d-2}({\bf s}) =0.
\ee
However, the conditions of B\'ezout's theorem are not satisfied because Proposition \ref{caustic} implies that the equations \eqref{maxdeg} have a common component consisting of the intersection $D_{n+1,0} =0 = D_{n,0}$ and  hence there are infinitely many solutions of \eqref{maxdeg}  in $[s_1:\dots :s_d]\in \mathbb P^{d-1}$ as long as $d\geq 4$.
On the other hand, as noted prior to the mentioned Proposition, this common component does not yield a meaningful solution to the Stieltjes--Fekete problem because $P_n(z)$ is the identically zero polynomial.   

In principle we should only count the solutions outside this locus, i.e. with $D_{n,0}\neq 0 $; a bit frustratingly, however, we cannot exclude that there are other common components, except that experiments suggest that this is not the case. This makes the application of general techniques of {\it Fulton's excess intersection formulas}\footnote{We thank Barbara Fantechi for pointing us towards this notion.} difficult to implement. A semi-heuristic argument (to be formalized in \cite{MasoeroFuturo}\footnote{We thank Davide Masoero for sharing this with us.}) would suggest that {\it generically}  the number of solutions should be the number of solutions of  the equation $k_1+k_2+\dots + k_d=n$ with $k_j$ non negative integers,  namely $(n+1)(n+2)\cdots (n+d-1)/(d-1)!$.

\begin{figure}
\begin{minipage}{0.5\textwidth}
\includegraphics[width=0.9\textwidth]{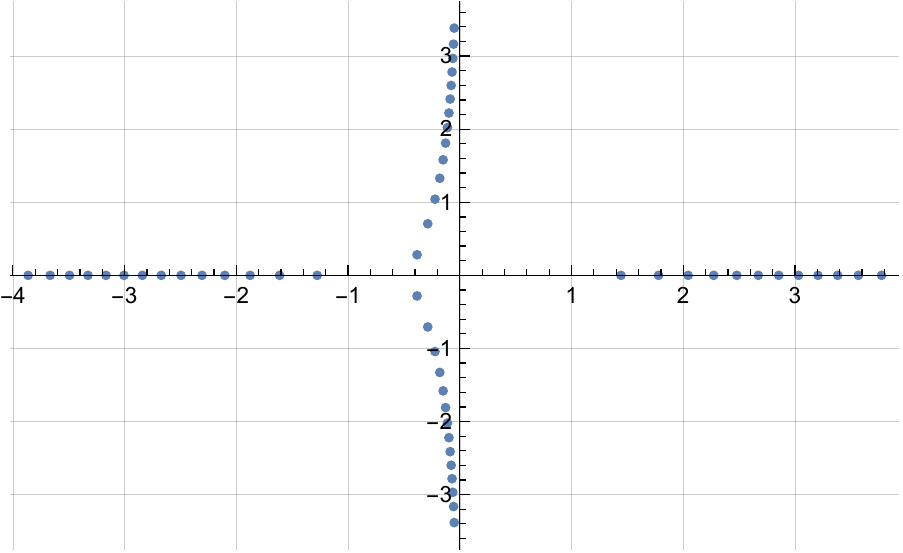}
\end{minipage}
\hfill
\begin{minipage}{0.44\textwidth}
\bea
&n=54\nn\\
&\frac {s_1}{s_3} \simeq -7.0607\,10^{-9} -2.6687\,10^{-7} i  \nn\\
&\frac{ s_2}{s_3} \simeq   0.998601\, -0.052877 i\nn
\\ \nn
&|\l_{\text{min}}(H_{n,0}) | \simeq  5.54711\, 10^{-5}
\\ \nn
&|\l_{\text{min}}(H_{n+1,0}) | \simeq  4.44482\, 10^{-59}
\\ \nn
&|\l_{\text{min}}(H_{n+1,1}) | \simeq  7.32318\, 10^{-59}
\eea
\end{minipage}
\\[14pt]
\begin{minipage}{0.5\textwidth}
\includegraphics[width=0.9\textwidth]{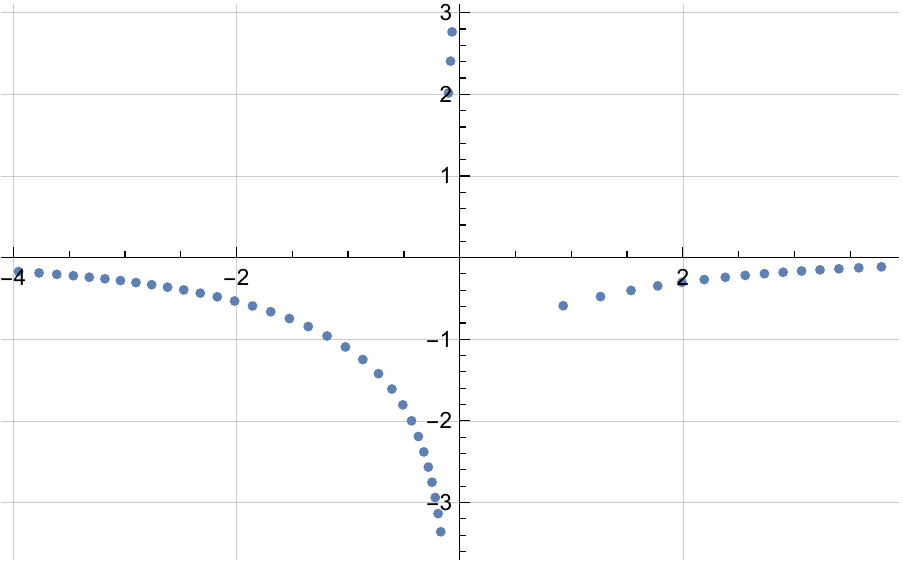}
\end{minipage}
\hfill
\begin{minipage}{0.44\textwidth}
\bea
&n=50\nn\\
&\frac {s_1}{s_3} \simeq -2.57799\,10^{-11} -5.50531\,10^{-12} i  \nn\\
&\frac{ s_2}{s_3} \simeq  -2.57799\,10^{-11} -5.50531\,10^{-12} i\nn
\\ \nn
&|\l_{\text{min}}(H_{n,0}) | \simeq 0.0874
\\ \nn
&|\l_{\text{min}}(H_{n+1,0}) | \simeq  4.98577\, 10^{-54}
\\ \nn
&|\l_{\text{min}}(H_{n+1,1}) | \simeq  5.19896\, 10^{-54}
\eea
\end{minipage}
\\[14pt]
\begin{minipage}{0.5\textwidth}
\includegraphics[width=0.9\textwidth]{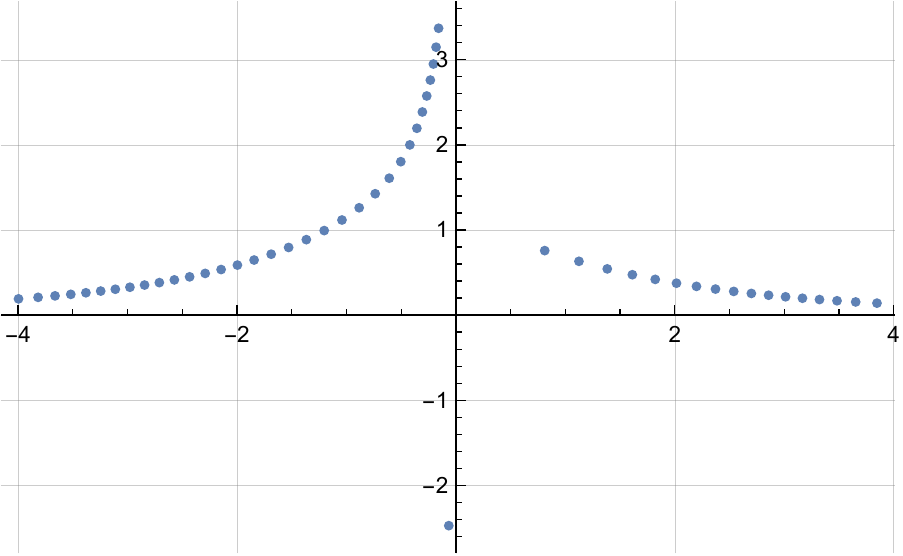}
\end{minipage}
\hfill
\begin{minipage}{0.44\textwidth}
\bea
&n=51\nn\\
&\frac {s_1}{s_2} \simeq -3.36574\,10^{-10} +1.02727\,10^{-10} i  \nn\\
&\frac{ s_3}{s_2} \simeq  -6.14402\,10^{-29} +2.32875\,10^{-28} i\nn
\\ \nn
&|\l_{\text{min}}(H_{n,0}) | \simeq 0.379591
\\ \nn
&|\l_{\text{min}}(H_{n+1,0}) | \simeq  6.74163\, 10^{-53}
\\ \nn
&|\l_{\text{min}}(H_{n+1,1}) | \simeq  4.86859\, 10^{-52}
\eea
\end{minipage}
\caption{Some solutions of the Bethe-Fekete equation \eqref{Fekete}; here $A=z^3-2z$ so that $\theta = -z^4/4 +z^2$. The plots indicate  particular solutions with $n$ points, obtained numerically. 
We then compute the values of $s_1,s_2,s_3$ using the formula \eqref{sj2}; only the ratios are relevant for the orthogonality and are indicated next to each panel figure. 
Using these values for the parameters $s_j$ we compute numerically the moments of the functional \eqref{Freudint}  and the determinants $D_{n+1,k}({\bf s})=\det H_{n+1,k}({\bf s})$,  $k=0,1$, in \eqref{Dnk} to verify that they are numerical zeros.  We then report the minimum of the absolute values of the eigenvalues of the matrices that appear in \eqref{Dnk}. Due to numerical approximations, the minimal eigenvalue is not exactly zero but since the computation was done with $50$ significant digits, they can be considered as vanishing. 
For comparison we report also  the  minimal eigenvalue for $H_{n,0}$ which is several orders of magnitude larger as expected since those matrices should not be singular.}
\label{Freudnumerics}
\end{figure}

{ 
\section{Numerical verifications}
Using Mathematica we can find numerical solutions of the equilibrium equations \eqref{Fekete} in some simple cases (it is required to use numerics in arbitrary precision) and then compute the resulting polynomial $P_n(z)$. 
One then uses the formulas \eqref{sj2}  for $s_j$ so as to obtain the moment functional.  We can then  compute the moments of the moment functional and plug them into the expressions for the  determinants $D_{n,k}$ \eqref{Dnk} to verify that they are (within numerical accuracy) indeed vanishing. 

We performed these numerical tests for the Freud case  with $A= z^3-2z$, $B=1$. 
The numerical solutions of \eqref{Fekete} are found using the ``FindRoot'' command, which will find the numerical solution of an algebraic system in a neighbourhood of an initially selected configuration by  Newton's iteration  method.

In this way we can test, with  different numerical solutions of \eqref{Fekete},  the whole setup. 
A few results of this numerical verification are reported in the gallery in Figure \ref{Freudnumerics}.
}

\section{Conclusions and further problems}
We conclude with a short discussion on the actual Fekete problem and some interesting connected problems that deserve further investigation.
\begin{enumerate}
\item[-]The core idea of this paper stems from \cite{BertoGravaHeredia} where we studied the numerical similarity between the exactly-integrable spectrum of the quartic anharmonic oscillator and the roots of the rational solutions of the second Painlev\'e\ transcendent. 

The setup in loc. cit. starts with a special case of the present result and requires the analysis of a particular type of maximally degenerate polynomials for a moment functional of type $(z^2-t,1)$; in this case the maximally degenerate polynomial is a $d=1$ degenerate one, namely, there is only one condition $D_{n,0}({\bf s})=0$ as in \eqref{Dnk}. 

However the crossing condition of the exactly-integrable spectrum of the corresponding anharmonic oscillator requires $d=2$ additional equations, namely 
\be
\int_{\gamma_j} P_n^2(z){ \rm e}^{\theta(z)} \d z=0, \ \ \ j = 1,2.
\label{QES}
\ee
Note that since $\mathcal M[P_n^2]=0$, there is, in fact, only one extra condition in \eqref{QES}. 
In the context of \cite{BertoGravaHeredia} this extra  condition was shown to be equivalent to the existence of a multiple eigenvalue for a suitable Sturm--Liouville boundary value problem and the extra condition was a constraint now on the values of $t$. 

Therefore, we can follow the same line of thought in this more  general setting and formulate the following problem:
for a maximally degenerate polynomial $P_n$ of type $(A,B)$, characterize the possible types $(A,B)$ for which the following  additional $d-1$ conditions hold:
\be
\int_{\gamma_j} P_n^2(z){ \rm e}^{\theta(z)} \d z=0, \ \ \ j = 1,2,\dots, d-1.
\label{QES+}
\ee
These ``exceptional types'' would be the direct analogues of the crossing conditions in the exactly integrable spectrum of the anharmonic oscillators that generalize the quartic example.

\item[-] The case $\deg A<\deg B$ and $B$ having only simple zeros corresponds to the classical Heine-Stiltjes polynomials and the analysis remains open. The contour $\Gamma$ is  a Pochhammer contours (group-commutators of the generators in the fundamental group of the plane minus the zeros of $B$) but there are  case distinctions according to whether the residues of $\theta'$ at these poles are positive integers, negative integers or neither which  complicate the description.

\item[-] The equilibrium equations \eqref{Fekete} imply a  hierarchy of similar equilibrium equation following ideas in \cite{Calogero78, Calogero77} used in the context of the  Stieltjes-Fekete problem associated to classical orthogonal polynomials. Briefly, the connection is as follows (in the simplest example of Hermite polynomials): the roots of the Hermite polynomials satisfy \eqref{Fekete} with $A(z)=z$ and $B\equiv 1$. The solutions can be thought of as the equilibrium (i.e. {\it stationary}) equations for the Hamiltonian 
\be
\frac 1 2 \sum_{j=1}^n (p_j^2 + x_j^2) - \sum_{j\neq \ell} \ln (x_j-x_\ell).
\ee
It was shown in loc. cit. that the same are {\it also} equilibrium solutions of the ``Calogero'' Hamiltonian;
\be
\frac 1 2 \sum_{j=1}^n (p_j^2 + x_j^2) + \sum_{j\neq \ell} \frac1{ (x_j-x_\ell)^2}.
\ee
In fact in \cite{Calogero78} they can be proved to be equilibrium equations of a hierarchy of Hamiltonians. 
{\color{black}
\item [-] The natural question arises regarding the asymptotic behaviour of the solutions to the Stieltjes-Bethe equations  for large $n$; the fact that these are the zeros of (degenerate) orthogonal polynomials opens the way to the use of the steepest descent method based on the Riemann--Hilbert analysis. However, the complication here is that the parameters $[s_1:\dots :s_d]$ depend in a rather  implicit way on $n$. An analysis in this direction was performed in the classical case of Heine--Stieltjes polynomials in \cite{Andrei2}, based on different approaches relying on potential theory. 
}
\item [-] Finally, the equations \eqref{Fekete} appear in the Bethe Ansatz, as already mentioned. 
For example \cite{HarnadWinternitz} they appear in the quantum separation of variables and hence our result implies that the solution of the Bethe Ansatz are particular cases of (semiclassical) orthogonal polynomials. This seems to be a significantly new result in a rather old subject and it is object of our future research.
\end{enumerate}

\paragraph{Acknowledgements.} 
M. B. would like to thank J. Harnad and F. Del Monte for useful discussions on the significance of the result in the context of the Bethe Ansatz.  
{\color{black} We also thank Davide Masoero for discussions about the number of solutions of the Stieltjes-Bethe  equations and Youn Miao  for pointing out several references on the Stieltjes-Bethe  equations.}

The work of M. B.  was supported in part by the Natural Sciences and Engineering Research Council of Canada (NSERC) grant RGPIN-2016-06660.
T.G. and E.C.H. acknowledge the funding from the European Union’s H2020 research and innovation programme under the Marie Sklodowska–Curie grant No. 778010 IPaDEGAN and the support of GNFM-INDAM group and the research project Mathematical Methods in NonLinear Physics (MMNLP), Gruppo 4-Fisica Teorica of INFN.   TG   acknowledges  the  hospitality and support from  Galileo Galilei Institute,
and from the scientific program on     "Randomness, Integrability, and Universality".

\appendix
\section*{Appendices}
\addcontentsline{toc}{section}{Appendices}
\renewcommand{\thesubsection}{\Alph{subsection}}
\renewcommand{\theequation}{\thesection.\arabic{equation}}
\section{One-parameter  family of solutions of Stieltjes-Bethe equations}
\label{family}
The following is an interesting example (see also the introduction in \cite{Boris11_1})  {because it provides an example where all zeros of $B$ are simple and also an example where the corresponding Stieltjes-Bethe equations admit a continuum of nontrivial solutions. At the same time this  also shows that in the case $\deg A\leq \deg B-1$ and with simple zeros of $B$, new phenomena occur (even under genericity assumptions) that are not present in the main theorems. }

Let $A =k B'$ with $k=1,2,\dots$,  so that ${\rm e}^{\theta} = \frac 1{B^{k+1}}$. Note that all zeros of $B$ must be simple for otherwise $A$ would not be relatively prime to $B$. 

The corresponding moment functionals are obtained by taking arbitrary linear combinations of small circles around the zeros of $B$; here we find  the first difference with the generic case, in that all of these $d+1$ moment functionals are linearly independent,  {which is a departure from the generic situation}. For example if $B=(z-\beta_1)(z-\beta_2)$ and $k=1$, then  the two linearly independent moment functionals $\mathcal M_{1,2}$ are
\be
\mathcal M_{\ell} [z^j] = \res{z=\beta_\ell } \frac {z^j\d z}{B^2(z)} = \frac {j \beta_\ell^ {j-1}}{ (\beta_1-\beta_2)^2} - 2 \frac {\beta_\ell^j}{(\beta_\ell - \beta_{3-\ell})^3},\ \ \ \ell = 1,2.
\ee
One can then verify that the first three moments $j=0,1,2$ form linearly independent vectors.

Now, the proofs of Theorems \ref{thm1} and \ref{thm2}  proceed as long as $n \geq \lfloor \frac {\deg B}2\rfloor +1$. 
For $n= k \deg B +1$ the Stieltjes--Bethe equations  \eqref{Fekete} have the immediate solution 
\be
P_n'(z) = B^k(z) \ \ \ P_n(z) = \int_0^z B^k(x) \d x + C
\ee
Here we see  that  the roots of $P_n$ provide a  one--parameter family of solutions parametrized by $C\in \C$ which, generically, consist of pairwise distinct points.  The polynomial $P$ is also the solution of 
\be
B P_n'' - kB' P_n' =0,
\ee
namely a generalized Lam\'e\ equation with trivial Van Vleck polynomial, according to the terminology in \cite{Boris11_1}.
Now the question is whether the polynomial $P_n$ is a maximally degenerate one; since $n = k\deg B+1$ we can retrace the proof of Theorem \ref{thm2} starting from 
\be
\label{24}
\int_\infty^z  \frac {B^k(w)\d w} {P_n^2(w)} = - \frac 1{P_n(z)} .
\ee
Let $\mathbb D_j$ be a small disk around $\beta_j$, \ \ $j=1,\dots, \deg B$. 
The function (see \eqref{intR})
\be
R_n(z):= \le\{
\begin{array}{lr}
\ds\frac{P_n(z)}{B(z)^{k+1}} \int_{\infty}^z\frac {B^k(w)\d w} {P_n^2(w)} =-\frac 1{B^{k+1}(z)} & z\not\in \cup_{j} \mathbb D_j\\[15pt]
\ds \frac{P_n(z)}{B(z)^{k+1}} \int_{\beta_j}^z\frac {B^k(w)\d w} {P_n^2(w)} =
  \frac {P_n(z)}{P_n(\beta_j)B(z)^{k+1}} - \frac 1{B(z)^{k+1}}   & z\in \mathbb D_j.
\end{array}
\ri.
\ee
The function satisfies clearly
\be
R_n(z_+)= R_n(z_-) + s_j \frac{P_n(z)}{B(z)^{k+1}},\ \ z\in \pa\mathbb D_j\  \ \ s_j:= \frac 1{P_n(\beta_j)}.
\ee
Thus 
\be
R_n(z)= \sum_{j=1}^{\deg B}  \frac 1{P_n(\beta_j)} \oint_{\pa\mathbb D_j}  \frac {P_n(w) \d w}{B^{k+1}(w)(w-z)2i\pi} 
\ee
and it follows that $P_n(z)$ is indeed a {\it family} of  maximally degenerate polynomials of degree $n = k \deg B+1$, whose zeros form a one-parameter family of solutions of the Stieltjes-Bethe equations \eqref{Fekete}, with $A = kB'$. 
In this case the moment functional is 
\be
\mathcal M[q(x)] =\sum_{j=1}^{\deg B}  \frac 1{P_n(\beta_j)} \oint_{\pa\mathbb D_j}  \frac {q(w) \d w}{B^{k+1}(w)2i\pi } .
\ee

{\color {blue}
\paragraph{The case $A=-kB'$, $k=1,2,\dots$, $B$ with simple roots.} 
In this case $\theta'=(k-1)\frac {B'}B$ and ${\rm e}^{\theta(z)} = B(z)^{k-1}$. According to our terminology the roots of $B$ are all {\it hard-edges} (for $k=1$) or {\it end-poles}  (for $k>1$), and contours $\gamma_1, \dots, \gamma_{d}$ (with $d=\deg B-1$) can be chosen for example as pairwise disjoint segments (except at endpoints) connecting one root of $B$ to the others. It follows from Cauchy's theorem that any other contour connecting two zeros  yields a functional that is the linear combination of the above ones. 

Also in this case the proof of the Theorem \ref{thm1} and \ref{thm2}  proceed as long as $n \geq \lfloor \frac {\deg B}2\rfloor +1$ without changes. 
}
{ 
\section{Lifting the primality condition on $A,B$}
\label{notprime}
The equilibrium  equations \eqref{Fekete} in Theorem \ref{thm1} are insensitive to the multiplication $A(z)\mapsto \wt A(z) = A(z) E(z)$ and $B(z)\mapsto \wt B(z)= B(z)E(z)$ for $E$ a polynomial. However this map does modify the definition of the moment functional for which $P_n(z) = \prod_{z=1}^N(z-z_j)$ is maximally degenerate orthogonal. On the face of it, this seems to be a conflict with our description, since the new moment functional of type $(\wt A, \wt B)$ has $\wt d = \max\{\deg \wt A, \deg \wt B-1\} = d+ \deg E$ and it is not clear that it admits the same maximally degenerate polynomial. 

We are now going to show that this is generically the case, namely,   the same  $P_n$ is maximally degenerate also for the modified moment functional provided that none of the  roots of $E$ coincides with one of $P_n$.

It suffices to consider the case $E(z) = (z-c)^K$; hence suppose that $\wt A(z) = A(z) (z-c)^K$ and $\wt B(z) = B(z)(z-c)^K$ (with $A,B$, relatively prime). 

The new moment functional $\wt \M$ is then related to the original $\M$ by 
\be
\wt \M[p(z)] = \sum_{j=1}^d s_j \int_{\gamma_j}p(z) \frac {{\rm e}^{\theta(z)}\d z}{(z-c)^K} + \sum_{\ell=1}^K (\ell-1)!\wt s_\ell \oint_{|z-c|=\epsilon}\frac {p(z)\d z}{(z-c)^\ell 2i\pi} 
\label{tildeM}
\ee
where $\gamma_j$ are the same contours defined earlier. Note that the new functional contains a linear combination of derivatives of the distributional delta function supported at $z=c$ (as a consequence of Cauchy's residue theorem).
To see that  \eqref{tildeM} holds, observe that the generating function of moments $\wt \Phi(\l):= \wt M[{\rm e}^{\l z}]$ satisfies the ODE (see \eqref{ODEgen})
\be
\bigg[\l B(\pa_\l) - A(\pa_\l) \bigg](\pa_\l-c)^K \wt \Phi(\l)=0
\ee
and at this point it is promptly verified that the general solution is 
\be
\wt \Phi(\l) = \sum_{j=1}^d s_j \int_{\gamma_j} \frac {{\rm e}^{\l z + \theta(z)}\d z}{(z-c)^K} + \sum_{\ell=1}^K (\ell-1)! \,\wt s_\ell  \oint_{|z-c|=\epsilon}\frac {{\rm e}^{\l z}\d z}{(z-c)^\ell 2i\pi} 
\ee
from which the statement follows. 

Suppose now that $P_n(z)$ is a maximally degenerate polynomial for $\M$, namely
\be
  0 = \M[P_n z^\ell] = \sum_{j=1}^d s_j \int_{\gamma_j}P_n(z)z^m  {{\rm e}^{\theta(z)}\d z}, \ \ m = 0,1,\dots, n-1, {\bf n, \dots, n+d-2},
\ee
for appropriate constants $s_1,\dots, s_d$. We claim that $P_n$ is also maximally degenerate with respect to $\wt M$ for  some additional constants $\wt s_1,\dots, \wt s_K$. 
Indeed, the new condition of maximal degeneracy now requires 
\be
0 = \wt \M[P_n z^\ell]=\sum_{j=1}^d s_j \int_{\gamma_j}P_n(z) z^m \frac {{\rm e}^{\theta(z)}\d z}{(z-c)^K} + \sum_{\ell=1}^K (\ell-1)!\,\wt s_\ell \oint_{|z-c|=\epsilon}\frac {P_n(z) z^m\d z}{(z-c)^\ell 2i\pi} \nn\\
m = 0, \dots, n-1, {\bf n, \dots, n+d-2+K}.
\ee
Since the sequence of polynomials $z^m$ and $(z-c)^m$ span the same flag of spaces, we can equivalently state the $\wt \M$-orthogonality of $P_n$ with respect to the sequence $(z-c)^m,  \ \ m = 0,\dots, n-2+d+K$. Of these $n-1+d+K$ conditions, the ones for $m = K,K+1,\dots, n+d-2+K$ are already satisfied automatically (since they reduce to the orthogonality of $P_n$ for the $\M$ moment functional). The $K$ conditions for $m=0,1,\dots, K-1$ fix  instead the constants $ \wt s_1,\dots, \wt s_K$: this defines a linear system for $\wt {\bf s}=(\wt s_1,\dots, \wt s_K)$ with a matrix of coefficients which is upper triangular and Toeplitz, of the form:
\be
 P_n(c)\1_{_K} +\sum_{\ell=1}^{K-1} P^{(\ell)}_n(c)\Lambda_{_K}^\ell,
\ee
with $\1_{_K}$ being the $K\times K$ identity matrix and $\Lambda_{_K}$ the upper shift matrix. If we choose $c$ generically (i.e. not equal to any of the zeros of $P_n$) then the system for the extra parameters $\wt {\bf s}$ is determined uniquely.
This shows that the new functional $\wt \M$ is also maximally degenerate on polynomials of order $n$ and that the maximally degenerate polynomial  of degree $n$ is the {\it same} as that for $\M$. 

If $c$ coincides with a zero of $P_n(c)$ in general one finds an inconsistent system because the orthogonality to $(z-c)^{K-1}$ yields the equation:
\be
\sum_{j=1}^d s_j \int_{\gamma_j}P_n(z)  \frac {{\rm e}^{\theta(z)}\d z}{(z-c)} = 0
\ee
which is never  satisfied because the expression is proportional by a non-vanishing function to the second linearly independent solution $G(z)$ (see \eqref{intR}) of the second-order ODE \eqref{SL}, and two independent solutions cannot vanish at the same point.
}

 \end{document}